# Towards regulator formulae for curves over number fields.

Rob de Jeu*

Department of Mathematical Sciences†,
University of Durham,
South Road,
Durham DH1 3LE,
United Kingdom

*In this paper we study the group $K_{2n}^{(n+1)}(F)$ where $F$ is the function field of a complete, smooth, geometrically irreducible curve $C$ over a number field, assuming the Beilinson–Soulé conjecture on weights. In particular, we compute the Beilinson regulator on a subgroup of $K_{2n}^{(n+1)}(F)$, using the complexes constructed in [dJ1]. We study the boundary map in the localization sequence for $n = 3$ (the case $n = 2$ was studied in [dJ2]). We combine our results with the results of Goncharov ([G1] and [G3], see also [G2]) in order to obtain a complete description of the image of the regulator map on $K_4^{(3)}(C)$ and $K_6^{(4)}(C)$, independent of any conjectures.*

## 1 Introduction

In [dJ2] we computed the regulator map and the boundary of a smooth, proper, geometrically irreducible curve $C$ over the number field $k$, on certain subgroups in $K_4^{(3)}(C)$ coming from the complexes defined by the author in [dJ1]. Combining this with the work of Goncharov [G1], one can obtain a complete description of the image under the regulator map of the whole of $K_4^{(3)}(C)$, see Theorem 5.4 and Corollary 5.5 below. One of the goals of the present paper is to extend the results of [dJ2] to subgroups of $K_{2n}^{(n+1)}(C)$ for all $n \geq 2$. However, although large parts can be done for arbitrary $n$, for $n \geq 4$ they can only be obtained assuming a standard conjecture about weights in algebraic $K$–theory. Similarly, although the approach taken here to computing the boundary map $\partial$ in the localization sequence

$$\cdots \longrightarrow K_{2n}^{(n+1)}(C) \longrightarrow K_{2n}^{(n+1)}(F) \xrightarrow{\partial} \coprod_{x \in C^{(1)}} K_{2n-1}^{(n)}(k(x)) \longrightarrow \cdots$$

on those subgroups probably works for any $n$, the combinatorics involved get rather complicated, so we carry it out completely only for $n = 3$. (The case $n = 2$ was done in [dJ2].) Putting our results together with the work of Goncharov we obtain a complete description of the image under the regulator map of $K_4^{(3)}(C)$ and $K_6^{(4)}(C)$, see Theorems 5.4 and 5.6 as well as Corollaries 5.5 and 5.8.

\* The author gratefully acknowledges support from HCM grant ERBCHBICT 941693.

† Current address: Dep. of Math., California Institute of Technology, Pasadena CA 91125, USA.

The paper is organized as follows. We review the construction of the complexes from [dJ1] in Section 2 below, state some of their properties, and take the opportunity to prove some loose ends needed in the rest of the paper. Section 3 contains the computation of the regulator map on the image of the cohomology groups in the $K$–theory of the field $F$, most easily described by pairing it with a holomorphic 1–form and integrating. Section 4 contains the computation of the boundary map on the image of the cohomology of the complexes, and is by far the longest. We give most of the proof for general $n$, but somewhere along the road the combinatorics simply become too complicated and we restrict ourselves to $n=3$. Finally, in Section 5 we relate our work with that of Goncharov, obtaining a complete combinatorical description of the image of $K_4^{(3)}(C)$ and $K_6^{(4)}(C)$ under the regulator map, independent of any conjectures. We conclude the section by indicating how such results could be obtained for higher $n$, but all this would depend rather heavily on conjectures in algebraic $K$–theory.

**Notation** The following notation will be fixed throughout the paper.

$k$ is a number field. $C$ is a smooth, geometrically irreducible, proper curve over $k$. $F = k(C)$ is the field of rational functions on $C$. $C^{(1)}$ will denote the set of points of $C$ of codimension one, i.e., the set of closed points of $C$.

In all sections except Section 2, $n$ is a fixed integer at least equal to two. For an Abelian group $A$, $A_\mathbb{Q} = A \otimes_\mathbb{Z} \mathbb{Q}$. $\mathbb{Q}(m) = (2\pi i)^m \mathbb{Q} \subset \mathbb{C}$ and similarly for $\mathbb{R}(m)$. In the decomposition $\mathbb{C} = \mathbb{R}(n-1) \oplus \mathbb{R}(n)$ we let $\pi_{n-1}$ denote the projection onto the $\mathbb{R}(n-1)$–part. If $S$ is a subset of a vector space $V$, we shall mean by $<S>$ (resp. $<S>_\mathbb{R}$) the $\mathbb{Q}$ (resp. $\mathbb{R}$) subspace spanned by the elements of $S$.

Throughout the paper, in integrals and cohomology groups, we write simply $C$ for $C_{\text{an}}$, which is the analytic manifold associated to $C \otimes_\mathbb{Q} \mathbb{C}$. Note that by our assuptions, this is a disjoint union of $[k:\mathbb{Q}]$ copies of a Riemann surface of genus the genus of $C$. Similarly, we shall write $H^1_{\text{dR}}(F)$ for $\varinjlim_{U \subset C} H^1_{\text{dR}}(C_{\text{an}})$ where the limit is over all Zariski open subsets of $C$.

We provide a little more motivation here for computing the regulators. Namely, the Beilinson conjectures for $C$ as above predict the following.

1) $K_{2n}^{(n+1)}(C)$ has $\mathbb{Q}$–dimension $r = \text{genus}(C)[k : \mathbb{Q}]$ for $n \geq 2$ (for $n = 1$, which we shall not study here, there is an additional integrality condition).

2) The Beilinson regulator induces an isomorphism

$$K_{2n}^{(n+1)}(C) \otimes_\mathbb{Q} \mathbb{R} \to H^2_\mathcal{D}(C; \mathbb{R}(n+1))^+ \cong H^1_{\text{dR}}(C; \mathbb{R}(n))^+.$$

3) If $\alpha_1, \ldots, \alpha_r$ is a $\mathbb{Q}$–basis of $K_{2n}^{(n+1)}(C)$, let $A$ be the matrix of writing $\text{reg}(\alpha_1), \ldots, \text{reg}(\alpha_r)$ with respect to a basis of $H^1_{\text{dR}}(C; \mathbb{Q}(n))^+ \cong \mathbb{Q}^r$. Note that $\det A$ is well defined as element of $\mathbb{R}^*/\mathbb{Q}^*$. Assume that the $L$–function $L(C, s)$ of $C$ can be analytically continued to the entire complex plane. Then $\det(A)/L^*(C, 1-n)$ is an element in $\mathbb{Q}^*$, where $L^*(C, n)$ is the first non–vanishing coefficient in the power series expansion of $L(C, s)$ around $s = n$.

Considering the role played by the image of the regulator in 2 and 3 above, it is important to have a good description of this image. In practice, it is often easier to determine the regulator by integrating $\text{reg}(\alpha) \wedge \overline{\omega}$ over $C$ for holomorphic 1–forms $\omega$ on $C$, and using the periods of the $\omega$'s, see Proposition 3.2 below.

## 2 Some preliminary results

This section contains a description of the complexes $\mathcal{M}_{(n)}^\bullet(F)$ and $\widetilde{\mathcal{M}}_{(n)}^\bullet(F)$, together with the maps $\varphi_{(n)}^p$ from their $H^p$ to $K_{2n-p}^{(n)}(F)$ under suitable assumptions, as they were constructed in [dJ1]. Apart from that, we also prove or state some results in this context that are useful for the rest of the paper.

We briefly recall the construction of the complexes $\mathcal{M}_{(n)}^\bullet(F)$ and $\widetilde{\mathcal{M}}_{(n)}^\bullet(F)$ in [dJ1, Section 3], where $F$ is a field of characteristic zero. Let $Y = \operatorname{Spec}(F)$, or more generally some reasonable regular scheme. Let $t$ be the standard affine coordinate on $\mathbb{P}^1$, and let $X_Y = \mathbb{P}_Y^1 \setminus \{t = 1\}$. In [dJ1] a formalism of "multi–relative" $K$–theory with weights is developed. To fix ideas, look at the exact sequence in relative $K$–theory

$$\ldots \to K_{m+1}^{(j)}(\{t=0,\infty\}) \to K_m^{(j)}(X_Y;\{t=0,\infty\}) \to K_m^{(j)}(X_Y) \to \ldots$$

One has $K_n^{(j)}(X_Y) \cong K_n^{(j)}(Y)$ by the homotopy property for $K$–theory of a reasonable regular scheme, and the map $K_n^{(j)}(X_Y) \to K_n^{(j)}(\{t=0,\infty\}) \cong K_n^{(j)}(Y)^{\oplus 2}$ is the diagonal map. From this, one gets isomorphisms $K_n^{(j)}(X_Y;\{t=0,\infty\}) \cong K_{n+1}^{(j)}(Y)$. (We shall apply this isomorphism only in case $Y$ is a Zariski open part of a smooth curve over a number field, the Spec of its function field, or the Spec of a number field, in which case all conditions are satisfied.) Iterating this idea one gets "multi-relative" $K$–theory, by taking relativity step by step. Let $t_i$ be the coordinate on the $i$–th copy of $X$ in $X^n$. Writing $\Box^n$ for $\{t_1 = 0, \infty\}, \ldots, \{t_n = 0, \infty\}$, we have a long exact sequence in relative $K$–theory

$$\cdots \to K_{m+1}^{(j)}(\{t_n=0,\infty\};\Box^{n-1}) \to K_m^{(j)}(X_Y^n;\Box^n) \to K_m^{(j)}(X_Y^n;\Box^{n-1}) \to \cdots$$

and as before it follows from the homotopy property for $K$–theory for some reasonable regular scheme $Y$ that $K_m^{(j)}(X_Y^n;\Box^n) \cong K_{m+1}^{(j)}(X_Y^{n-1};\Box^{n-1})$ for $m \geq 0$. Repeating this, we get $K_m^{(j)}(X_Y^n;\Box^n) \cong K_{m+n}^{(j)}(Y)$ for $m \geq 0$. Note that there is no obvious choice of this isomorphism, which will result in statements up to sign below.

Let $Y = \operatorname{Spec}(F)$, but drop $Y$ from the notation. Let $U \subset F^* \setminus \{1\}$ be finite. Write $X_{\operatorname{loc}}^k = X^k \setminus \{t_i = u_j, \ u_j \in U, \ i = 1, \ldots, k\}$. One has a fourth quadrant spectral sequence

$$E_1^{p,q} = \coprod K_{-p-q}^{(n-p)}(X_{\operatorname{loc}}^{n-1-p};\Box^{n-1-p}) \Rightarrow K_{-p-q}^{(n)}(X^{n-1};\Box^{n-1}) \cong K_{-p-q+n-1}^{(n)}(F) \quad (2.1)$$

which looks like

$$\begin{array}{lll}
K_{-q-1}^{(n)}(X_{\operatorname{loc}}^{n-1};\Box^{n-1}) & \coprod K_{-q-2}^{(n-1)}(X_{\operatorname{loc}}^{n-2};\Box^{n-2}) & \coprod K_{-q-3}^{(n-2)}(X_{\operatorname{loc}}^{n-3};\Box^{n-3}) \\
K_{-q}^{(n)}(X_{\operatorname{loc}}^{n-1};\Box^{n-1}) & \coprod K_{-q-1}^{(n-1)}(X_{\operatorname{loc}}^{n-2};\Box^{n-2}) & \coprod K_{-q-2}^{(n-2)}(X_{\operatorname{loc}}^{n-3};\Box^{n-3}) \\
K_{-q+1}^{(n)}(X_{\operatorname{loc}}^{n-1};\Box^{n-1}) & \coprod K_{-q}^{(n-1)}(X_{\operatorname{loc}}^{n-2};\Box^{n-2}) & \coprod K_{-q-1}^{(n-2)}(X_{\operatorname{loc}}^{n-3};\Box^{n-3})
\end{array}$$

Here the coproduct for $X_{\operatorname{loc}}^{n-1-p}$ corresponds to the codimension $p$ hyperplanes given by $p$ equations of type $t_i = u_i$, $u_i \in U$. If $K_m^{(j)}(Y) = 0$ for $2j \leq m$, $m > 0$, all the terms below the row with $q = -n$ vanish, [dJ1, page 221]. Hence if we view this lowest row with the differential of the spectral sequence as a cohomological complex (depending on $U$)

$$C_{(n)}^\bullet: \ K_n^{(n)}(X_{\operatorname{loc}}^{n-1};\Box^{n-1}) \to \coprod K_{n-1}^{(n-1)}(X_{\operatorname{loc}}^{n-2};\Box^{n-2}) \to \coprod K_{n-2}^{(n-2)}(X_{\operatorname{loc}}^{n-3};\Box^{n-3}) \to \ldots$$

in degrees 1 through $n$, we get a map

$$H^p(C^\bullet_{(n)}) \to K^{(n)}_{n-p+1}(X^{n-1}; \square^{n-1}) \cong K^{(n)}_{2n-p}(F).$$

This procedure works more generally for $Y$ a reasonable regular scheme, and $U \subset \Gamma(Y, \mathcal{O}^*) \setminus \{\infty\}$ such that for all $u_k$ and $u_l$ in $U$, $u_k - u_l$ and $u_k - 1$ are invertible on $Y$ if they are not identically zero. Let $G = \operatorname{Spec}(\mathbb{Q}[S, S^{-1}, (1-S)^{-1}])$. From the localization sequence

$$\cdots \to K^{(j-1)}_n(\operatorname{Spec}(\mathbb{Q}))^\oplus \to K^{(j)}_n(\mathbb{A}^1_\mathbb{Q}) \to K^{(j)}_n(G) \to K^{(j-1)}_{n-1}(\operatorname{Spec}(\mathbb{Q}))^{\oplus 2} \to \cdots$$

and the facts that $K^{(j)}_n(\mathbb{A}^1_\mathbb{Q}) \cong K^{(j)}_n(\mathbb{Q})$ and $K^{(j)}_n(\mathbb{Q}) = 0$ unless $n = 2j-1$ for $j \geq 1$, one gets that the conditions about weights above are satisfied for $G$. One can use the spectral sequence above, with $G$ instead of $Y$, and $U = \{S\}$, to construct elements $[S]_n \in K^{(n)}_n(X^{n-1}_{G,\mathrm{loc}}; \square^{n-1})$ for $n \geq 2$, satisfying $\mathrm{d}[S]_n = \sum_{j=1}^{n-1}(-1)^j [S]_{n-1|t_j=S}$, where we put $[S]_1 = 1 - S$. With some more care, one sees that actually $[S]_n \in K^{(n)}_n(X^{n-1}_{\mathbb{G}_m,\mathrm{loc}}; \square^{n-1})$. Any $u \in F^* \setminus \{1\}$, or more generally any $u \in \Gamma(Y, \mathcal{O}^*)$ such that $1 - u$ is also invertible on $Y$, yields a map $Y \to G$, and hence yields an element $[u]_n \in K^{(n)}_n(X^{n-1}_{Y,\mathrm{loc}}; \square^{n-1})$ by pulling back, with boundary $\mathrm{d}[u]_n = \sum_{j=1}^{n-1}(-1)^j [u]_{n-1|t_j=u}$ in $C^\bullet_{(n)}$.

We now return to the case $Y = \operatorname{Spec}(F)$, $U \subset F^* \setminus \{1\}$ finite, $X^n_{\mathrm{loc}}$ as before. Write $K_{(p)}$ for $K^{(p)}_p(X^{p-1}_{Y,\mathrm{loc}}; \square^{p-1})$. For the construction of $\widetilde{\mathcal{M}}^\bullet_{(n)}((F))$ one starts with the complex $C^\bullet_{(n)}$ (starting in degree 1)

$$K_{(n)} \to \left(\coprod_{U'} K_{(n-1)}\right)^{\oplus \binom{n-1}{1}} \to \left(\coprod_{U'^2} K_{(n-2)}\right)^{\oplus \binom{n-1}{2}} \to \cdots \to \left(\coprod_{U'^{n-1}} K_{(1)}\right)^{\oplus \binom{n-1}{n-1}}.$$

The $\oplus \binom{n-1}{p}$ here corresponds to the directions number of ways of putting $p$ of the coordinate $t_j$ to a constant in $U'$. For any $u \in U$, we have $[u]_n \in K^{(n)}_n(X^{n-1}_{\mathrm{loc}}; \square^{n-1})$. The element $[u]_2$ has boundary $(1-u)^{-1}_{|t=u}$, and for $n \geq 3$ $[u]_n$ has boundary $\sum_{j=1}^{n-1}(-1)^j [u]_{n-1|t_j=u}$. Moreover, $C^\bullet_{(n)}$ carries an action of $S_{n-1}$ by permuting the coordinates, and $[u]_n$ is in fact in the alternating part for this action. Let

$$(1+I)^* = K^{(1)}_1(X_{\mathrm{loc}}; \square) = \left\{ F(t) = \prod_i (t-x_i)^{n_i}(t-1)^{-n_i} \,\big|\, x_i \in U, \, \prod_i x_i^{n_i} = 1 \right\}_\mathbb{Q}.$$

There are $m-1$ cup products

$$(1+I)^* \cup K^{(m-1)}_{m-1}(X^{m-2}_{\mathrm{loc}}; \square^{m-2}) \to K^{(m)}_m(X^{m-1}_{\mathrm{loc}}; \square^{m-1})$$

depending on which of the coordinates on $X^{m-1}_{\mathrm{loc}}$ we use for the $(1+I)^*$–factor. We let $(1+I)^* \tilde{\cup} K^{(m-1)}_{m-1}(X^{m-2}_{\mathrm{loc}}; \square^{m-2})$ denote the span of the images of all possibilities. Define

$$\mathrm{symb}_2 = <[u]_2> + (1+I)^* \cup F^*_\mathbb{Q} \subset K^{(2)}_2(X_{\mathrm{loc}}; \square)$$

and for $n \geq 3$

$$\mathrm{symb}_n = <[u]_n> + (1+I)^*\tilde{\cup}\,\mathrm{symb}_{n-1} \subset K_n^{(n)}(X_{\mathrm{loc}}^{n-1}; \square^{n-1}).$$

We get a subcomplex $C^\bullet_{(n),\log}$ of $C^\bullet_{(n)}$,

$$\mathrm{symb}_n \to \left(\coprod_{U'}\mathrm{symb}_{n-1}\right)^{\oplus\binom{n-1}{1}} \to \left(\coprod_{U'^2}\mathrm{symb}_{n-2}\right)^{\oplus\binom{n-1}{2}} \to \cdots \to \left(\coprod_{U'^{n-1}} F^*_\mathbb{Q}\right)^{\oplus\binom{n-1}{n-1}}.$$

The subcomplex $J^\bullet_{(n)}$ of $C^\bullet_{(n),\log}$ given by

$$(1+I)^*\tilde{\cup}\,\mathrm{symb}_{n-1} \to \mathrm{d}(\ldots) + \left(\coprod (1+I)^*\tilde{\cup}\,\mathrm{symb}_{n-1}\right)^{\oplus\binom{n-1}{1}} \to \mathrm{d}(\ldots) + \ldots \to \cdots$$

$$\cdots \to \mathrm{d}(\ldots) + \left(\coprod (1+I)^*\tilde{\cup}\,\mathrm{symb}_2\right)^{\oplus\binom{n-1}{n-2}} \to \mathrm{d}(\ldots)$$

is acyclic, and we can form the quotient complex $C^\bullet_{(n),\log}/J^\bullet_{(n)}$. Because $S_{n-1}$ acts on $C^\bullet_{(n),\log}$ and $J^\bullet_{(n)}$ is stable under the action, we can take the alternating part of this quotient complex, and we get the complex

$$\mathcal{M}^\bullet_{(n)}(F): \quad M_{(n)} \to M_{(n-1)} \otimes <U> \to \cdots \to M_{(2)} \otimes \bigwedge^{n-2} <U> \to F^*_\mathbb{Q} \otimes \bigwedge^{n-1} <U>.$$

where $<U>$ is the (multiplicative) subspace of $F^*_\mathbb{Q}$ spanned by $U$, and

$$M_{(k)} = \frac{\mathrm{symb}_k}{(1+I)^*\tilde{\cup}\,\mathrm{symb}_{k-1}}.$$

(In [dJ1] and [dJ2] we wrote the factors in the tensor product the other way round. We change this notation here to conform with the notation used by Goncharov). Finally, by taking direct limits over $U$ we get the complex

$$\mathcal{M}^\bullet_{(n)}(F): \quad M_{(n)}(F) \to M_{(n-1)}(F) \otimes F^*_\mathbb{Q} \to \cdots \to M_{(2)}(F) \otimes \bigwedge^{n-2} F^*_\mathbb{Q} \to F^*_\mathbb{Q} \otimes \bigwedge^{n-1} F^*_\mathbb{Q}.$$

So $M_{(k)}(F)$ is generated by symbols $[f]_k$ with $f \in F^* \setminus \{1\}$, and the differential is given by

$$\mathrm{d}\left([f]_k \otimes g_1 \wedge \cdots \wedge g_{n-k}\right) = [f]_{k-1} \otimes f \wedge g_1 \wedge \cdots \wedge g_{n-k}$$

if $k \geq 3$ and

$$(1-f) \otimes f \wedge g_1 \wedge \cdots \wedge g_{n-k}$$

if $k = 2$. Note that the symbol $[1]_k$ also exists (with the relation $[1]_k = 2^{k-1}([1]_k + [-1]_k)$ as proved in [dJ1, Proposition 6.1]), and in particular $\mathrm{d}[1]_k = 0$ for all $k \geq 2$. By construction, if the Beilinson–Soulé conjecture holds for $F$, there are maps

$$\varphi^p_{(n)} : H^p(\mathcal{M}^\bullet_{(n)}(F)) \to K^{(p)}_{2n-p}(F)$$

as the composition of

$$H^p(\mathcal{M}^\bullet_{(n)}(F))\stackrel{\sim}{\leftarrow} H^p(C^\bullet_{(n),\log}(F)^{\mathrm{alt}}) \to H^p(C^\bullet_{(n),\log}(F)) \to H^p(C^\bullet_{(n)}(F)) \to K^{(p)}_{2n-p}(F).$$

Finally, the complex $\widetilde{\mathcal{M}}^\bullet_{(n)}(F)$ is obtained by quotienting out the complex $\mathcal{M}^\bullet_{(n)}(F)$ by the subcomplex

$$<[u]_n + (-1)^n[1/u]_n> \to <[u]_{n-1} + (-1)^{n-1}[1/u]_{n-1}> \otimes F^*_{\mathbb{Q}} \to \cdots \quad (2.2)$$

$$\cdots \to <[u]_2 + (-1)^n[1/u]_2> \otimes \bigwedge^{n-2} F^*_{\mathbb{Q}} \to \mathrm{d}(\ldots).$$

We get the complex $\widetilde{\mathcal{M}}^\bullet_{(n)}(F)$

$$\widetilde{M}_{(n)}(F) \to \widetilde{M}_{(n-1)}(F) \otimes F^*_{\mathbb{Q}} \to \cdots \to \widetilde{M}_{(2)}(F) \otimes \bigwedge^{n-2} F^*_{\mathbb{Q}} \to \bigwedge^n F^*_{\mathbb{Q}}.$$

where $\widetilde{M}_{(k)}(F) = M_{(k)}(F)/<[u]_k+(-1)^k[1/u]_k>$. The subcomplex (2.2) is acyclic in degrees $n-1$ and $n$ ([dJ1, Remark 3.23]) and is acyclic everywhere if the Beilinson–Soulé conjecture is true (not just for $K_*(F)$ but for more schemes, see [dJ1, Proposition 3.20]). Note that now the differential at the $(n-1)$–th place is given by

$$\mathrm{d}\left([f]_2 \otimes g_1 \wedge \cdots \wedge g_n\right) = (1-f) \wedge f \wedge g_1 \wedge \cdots \wedge g_n$$

with the other differentials unchanged. If the Beilinson–Soulé conjecture holds more generally, we therefore get a map

$$\tilde\varphi^p_{(n)} : H^p(\widetilde{\mathcal{M}}^\bullet_{(n)}(F)) \to K^{(p)}_{2n-p}(F)$$

as the composition of

$$H^p(\widetilde{\mathcal{M}}^\bullet_{(n)}(F))\stackrel{\sim}{\leftarrow} H^p(\mathcal{M}^\bullet_{(n)}(F))\stackrel{\sim}{\leftarrow} H^p(C^\bullet_{(n),\log}(F)^{\mathrm{alt}}) \to H^p(C^\bullet_{(n)}(F)) \to K^{(p)}_{2n-p}(F)$$

Here the leftmost arrow is an isomorphism if the Beilinson-Soulé conjecture is true in general, and the rightmost arrow exists if the Beilinson-Soulé conjecture is true for the $K$–theory of $F$. By construction, all arrows from left to right are injective for $p = 1$, if they exist.

The reader may check that, if we assume the Beilinson–Soulé conjecture in general, for any element $\alpha$ in $H^p(\mathcal{M}^\bullet_{(n)}(F))$ (resp. $H^2(\widetilde{\mathcal{M}}^\bullet_{(n)}(F))$), $\varphi^p_{(n)}(\alpha)$ (resp. $\tilde\varphi^p_{(n)}(\alpha)$) naturally lives in $K^{(n)}_{2n-p}(U)$ for $U$ some Zariski open subset of a "reasonable" regular scheme $Y$ with function field $F$. This is because the lift of such an element will involve only finitely many elements in $F$, and the spectral sequence (2.1) will involve only finitely many $t_i = u_j$'s. But then this spectral sequence exists for a suitable open part of $Y$ as well, by leaving out the closed part where $u_i = u_j$ for all $i,j$ such that $u_i \neq u_j$.

In this paper we shall be mainly interested in the case $p = 2$ and $n = 4$, i.e., the target is $K^{(4)}_6(F)$. The leftmost arrow here is a surjection without any assumptions because of the acyclicity of the complex (2.2) in degree 3. The rightmost arrow exists to a quotient $K^{(4)}_6(F)/N$, which is as follows. In the spectral sequence (2.1) all higher differentials leaving $E_2^{1,-4}$ are zero, as they land in $K_2^{(1)}(F)$'s or outside the range of the spectral sequence. So $E_2^{1,-4} = E_\infty^{-1,4}$ and we get a map $H^2(C_{(4)}^\bullet(F)) = E_\infty^{-1,4}$, a subquotient of $K_6^{(4)}(F)$. In order to determine this more precisely, note that we have a long exact localization sequence

$$\cdots \to \coprod K_3^{(1)}(F) \to K_3^{(2)}(X;\square) \to K_3^{(2)}(X_{\text{loc}};\square) \to \coprod K_2^{(1)}(F) \to \cdots.$$

As $K_3^{(1)}(F)$ and $K_2^{(1)}(F)$ are both zero, we get

$$K_3^{(2)}(X_{\text{loc}};\square) \cong K_3^{(2)}(X;\square) \cong K_4^{(2)}(F).$$

Therefore $E_\infty^{-1,4} \subset K_6^{(4)}(F)/N$ with $N$ generated by $K_4^{(2)}(F) \cup K_2^{(2)}(F)$, and we get a map

$$H^2(\mathcal{M}_{(4)}^\bullet(F)) \to E_\infty^{1,-4} \to K_6^{(4)}(F)/N$$

which does not depend on any assumptions.

In Proposition 4.1 below, we shall introduce maps $\delta = \coprod \delta_x$ with

$$\delta_x : \tilde{\mathcal{M}}_{(n)}^\bullet(F) \to \widetilde{\mathcal{M}}_{(n-1)}^\bullet(k(x))[1]$$

given by

$$\delta_x([f]_{n-k} \otimes g_1 \wedge \ldots \wedge g_k) = \text{sp}_{n-k,x}([f]_{n-k}) \otimes \partial_{k,x}(g_1 \wedge \ldots \wedge g_k)$$

in degrees 1 through $n-1$, and by $-\partial_{n,x}$ in degree $n$, where $\text{sp}_{n-k,x}([f]_{n-k}) = [f(x)]_{n-k}$ if $f(x) \neq 0$ or $\infty$, 0 otherwise, and $\partial_{k,x}$ the unique map from $\bigwedge^k F_\mathbb{Q}^*$ to $\bigwedge^{k-1} k(x)_\mathbb{Q}^*$ determined by

$$\pi_x \wedge u_1 \wedge \ldots \wedge u_{k-1} \mapsto u_1(x) \wedge \cdots \wedge u_{k-1}(x)$$
$$u_1 \wedge \ldots \wedge u_k \mapsto 0$$

if all $u_i$ are units at $x$ and $\pi_x$ is a uniformizer at $x$. This map obviously gives rise to a map $\widetilde{\mathcal{M}}_{(n)}^\bullet(F) \to \widetilde{\mathcal{M}}_{(n-1)}^\bullet(k(x))[1]$ by composition with the natural projection $\mathcal{M}_{(n)}^\bullet(F) \to \widetilde{\mathcal{M}}_{(n)}^\bullet(F)$. Following Goncharov ([G2, § 6]), we introduce the complexes $\mathcal{M}_{(n+1)}^\bullet(C)$ and $\widetilde{\mathcal{M}}_{(n+1)}^\bullet(C)$, defined to be the total complexes of

$$\begin{array}{ccccccccc}
0 & \to & M_{(n+1)}(F) & \xrightarrow{d} & M_{(n)}(F) \otimes F_\mathbb{Q}^* & \xrightarrow{d} & M_{(n-1)}(F) \otimes \bigwedge^2 F_\mathbb{Q}^* & \xrightarrow{d} & \cdots \\
& & \downarrow & & \delta \downarrow & & \delta \downarrow & & \\
& & 0 & \to & \coprod \widetilde{M}_{(n)}(k(x)) & \xrightarrow{d} & \coprod \widetilde{M}_{(n-2)}(k(x)) \otimes k(x)_\mathbb{Q}^* & \xrightarrow{d} & \cdots
\end{array}$$

and

$$\begin{array}{ccccccccc}
0 & \to & \widetilde{M}_{(n+1)}(F) & \xrightarrow{d} & \widetilde{M}_{(n)}(F) \otimes F_\mathbb{Q}^* & \xrightarrow{d} & \widetilde{M}_{(n-1)}(F) \otimes \bigwedge^2 F_\mathbb{Q}^* & \xrightarrow{d} & \cdots \\
& & \downarrow & & \delta \downarrow & & \delta \downarrow & & \\
& & 0 & \to & \coprod \widetilde{M}_{(n)}(k(x)) & \xrightarrow{d} & \coprod \widetilde{M}_{(n-2)}(k(x)) \otimes k(x)_\mathbb{Q}^* & \xrightarrow{d} & \cdots
\end{array}$$

where both coboundaries have degree 1 and the total complexes are cohomological complexes with $\mathcal{M}_{(n+1)}^\bullet(F)$ and $\widetilde{\mathcal{M}}_{(n+1)}^\bullet(F)$ in degree 1. There is an obvious inclusion of $H^2(\mathcal{M}_{(n+1)}^\bullet(C))$ into $H^2(\mathcal{M}_{(n+1)}^\bullet(F))$, and similarly we have an inclusion of $H^2(\widetilde{\mathcal{M}}_{(n+1)}^\bullet(C))$

into $H^2(\widetilde{\mathcal{M}}^\bullet_{(n+1)}(F))$, so that the maps $\varphi^2_{(n+1)}$ resp. $\tilde\varphi^2_{(n+1)}$ obviously extend to maps on the cohomology of those complexes.

In [dJ1] regulator maps

$$K_p^{(q)}(X^n_{Y,\mathrm{loc}};\square^n) \to H_\mathcal{D}^{2q-p}(X^n_{Y,\mathrm{loc}};\square^n;\mathbb{R}(q))$$

to relative Deligne cohomology were defined. We recall that

$$H_\mathcal{D}^n(X;E;\mathbb{R}(q)) \cong \frac{\left\{\begin{array}{l}(\omega_n,s_n) \text{ with } \omega_n \in F^q(D)^n,\ s_n \in j_*S_X^{n-1}(q-1) \text{ such}\\ \text{that } \omega_{n|E}\equiv 0, s_{n|E}\equiv 0 \text{ and } \mathrm{d}s_n = \pi_{q-1}\omega_n\end{array}\right\}}{\left\{\begin{array}{l}(\mathrm{d}\omega_{n-1}, \pi_{q-1}\omega_{n-1}-\mathrm{d}s_{n-1}) \text{ with } \omega_{n-1}\in F^q(D)^{n-1},\\ s_{n-1}\in j_*S_X^{n-2}(q-1) \text{ such that } \omega_{n-1|E}\equiv 0, s_{n-1|E}\equiv 0\end{array}\right\}}.$$

(See [dJ1, p. 218].) Here the notation means the following. We write $X$ etc. for the underlying topological complex manifold consisting of the closed points of $X\times_{\mathrm{Spec}\,(\mathbb{Q})}\mathrm{Spec}\,(\mathbb{C})$. $\overline{X}$ is a compactification of $X$ with complement $D$ such that $D$ and $D\cup E$ are a system of divisors with normal crossings. $j$ is the imbedding of $X$ into $\overline{X}$. $S_X^\bullet(q)$ is the complex of $\mathbb{R}(q)$–valued $C^\infty$–forms on $X$, $F^q(D)^\bullet$ the complex of $\mathbb{C}$–valued $C^\infty$–forms on $\overline{X}$ of type $(p,r)$ with $p\geq q$ and with logarithmic poles along $D$. (So locally on $\overline{U}\subset \overline{X}$ an element in $F^q(D)^n$ is a sum of elements of the form $\varphi\wedge\psi$ with $\varphi\in\Omega^\bullet_{\overline{U}}(D\cap\overline{U})$ of degree $p\geq q$, and $\psi\in C^{0,n-p}(\overline{U})$.) Note that if $q>\dim X$, we get a natural isomorphism $H_\mathcal{D}^n(\mathbb{R}(q))\cong H_{\mathrm{dR}}^{n-1}(\mathbb{R}(q-1))$.

The regulator lands in the invariant (or plus) part of Deligne cohomology with respect to the involution given by the combined action of complex conjugation on the underlying topological space (through the action on $\mathbb{C}$ in $X\times_{\mathrm{Spec}\,(\mathbb{Q})}\mathrm{Spec}\,(\mathbb{C})$) and on the coefficients $\mathbb{R}(q)\subset\mathbb{C}$. This involution acts similarly on $H_{\mathrm{dR}}^{n-1}(\mathbb{R}(q-1))$, and the plus–space in Deligne cohomology is isomorphic to the plus–space in $H_{\mathrm{dR}}$ if $q>\dim X$.

The regulator of a cup product in $K$–theory is given by the cup products of the regulators, (see [dJ1, (22) and (40)], but (40) is flawed by typographical errors). For $(\omega_p,s_p)$ in $H_\mathcal{D}^p(\mathbb{R}(k))$, $(\omega_q,s_q)$ in $H_\mathcal{D}^q(\mathbb{R}(l))$, we get that in $H_\mathcal{D}^{p+q}(\mathbb{R}(k+l))$

$$(\omega_p,s_p)\cup(\omega_q,s_q) = (\omega_p\wedge\omega_q, s_p\wedge\pi_l\omega_q + (-1)^p(\pi_k\omega_p)\wedge s_q) \qquad (2.3)$$

As for the regulator of $[S]_n$, it is given by $(\omega_n,\varepsilon_n)$, with

$$\omega_n = (-1)^{n-1}\mathrm{d}\log\frac{t_1-S}{t_1-1}\wedge\cdots\wedge\mathrm{d}\log\frac{t_{n-1}-S}{t_{n-1}-1}\wedge\mathrm{d}\log(1-S).$$

Here $\varepsilon_n$ is an $\mathbb{R}(n-1)$–valued $(n-1)$–form such that $\mathrm{d}\varepsilon_n = \pi_{n-1}\omega_n$. (Unfortunately the signs in equation (41) in [dJ1] were wrong, so the formula for $\omega_{p+1}$ on page 237 needs a sign $(-1)^p$. this does not change the results of the paper, as it only introduces a similar sign in [dJ1, Proposition 4.1], which was stated up to sign anyway.)

Finally, we have to introduce some polylogarithm functions and state their relations with the present constructions.

Let $\mathrm{Li}_k(z) = \sum_{m=1}^\infty z^m/m^k$ for $k\geq 1$ and $z\in\mathbb{C}$, $|z|<1$. Then $\mathrm{Li}_1(z) = -\mathrm{Log}(1-z)$ and $\mathrm{dLi}_{k+1}(z) = \mathrm{Li}_k(z)\mathrm{d}\log z$ for $k\geq 1$. The functions $\mathrm{Li}_k$ can be continued to multi–valued holomorphic functions on $\mathbb{P}^1_\mathbb{C}\setminus\{0,1,\infty\}$. Let the Bernoulli numbers $B_k$ be defined by

$$\frac{x}{e^x-1} = \sum_{k=0}^\infty \frac{B_k}{k!}x^k$$

It is well known that the functions (called $P_n$ resp. $P_{n,\text{Zag}}$ in [dJ1] and [dJ2])

$$P_{n,\text{Zag}}(z) = \pi_{n-1}\left(\sum_{k=0}^{n-1} \frac{(-\log|z|)^k}{k!}\text{Li}_{n-k}(z)\right)$$

and

$$P_n^{\text{mod}}(z) = \pi_{n-1}\left(\sum_{k=0}^{n-1} \frac{2^k B_k}{k!}\log^k|z|\text{Li}_{n-k}(z)\right)$$

extend to single valued functions on $\mathbb{P}_{\mathbb{C}}^1 \setminus \{0,1,\infty\}$ with values in $\mathbb{R}(n-1)$, see [Za]. The functions $P_n^{\text{mod}}$ satisfy the functional equations $P_n^{\text{mod}}(z) + (-1)^n P_n^{\text{mod}}(z^{-1}) = 0$, and extend to continuous functions on $\mathbb{P}_{\mathbb{C}}^1$ with $P_n^{\text{mod}}(0) = P_n^{\text{mod}}(\infty) = 0$.

We have the following relations between the functions $P_{n,\text{Zag}}$:

**Lemma 2.1**

$$dP_{n,\text{Zag}}(z) = P_{n-1,\text{Zag}}(z)d\,i\arg z + (-1)^n \frac{\log^{n-1}|z|}{(n-1)!}\pi_{n-1}d\log(1-z)$$

**Proof**

$$\begin{aligned}dP_{n,\text{Zag}}(z) =&\pi_{n-1}\left(\sum_{k=0}^{n-2}\frac{(-\log|z|)^k}{k!}\text{Li}_{n-1-k}(z)d\log z\right)\\ &+ \pi_{n-1}\left(-\frac{(-\log|z|)^{n-1}}{(n-1)!}d\log(1-z)\right)\\ &- \pi_{n-1}\left(\sum_{k=1}^{n-1}\frac{(-\log|z|)^{k-1}}{(k-1)!}\text{Li}_{n-k}(z)d\log|z|\right)\\ =&\pi_{n-1}\left(\sum_{k=0}^{n-2}\frac{(-\log|z|)^k}{k!}\text{Li}_{n-1-k}(z)d\,i\arg z\right)\\ &+ (-1)^n\frac{\log^{n-1}|z|}{(n-1)!}\pi_{n-1}d\log(1-z)\end{aligned}$$

As in [Za, §7] one checks that we have the relations

$$P_{n,\text{Zag}}(z) = \sum_{0\leq 2j < m} \frac{\log^{2j}|z|}{(2j+1)!} P_{m-2j}^{\text{mod}}(z).$$

**Lemma 2.2** Let $C$ be a complete, smooth, irreducible curve over $\mathbb{C}$ with function field $F = \mathbb{C}(C)$. If $f_1, \ldots, f_k$ are elements of $F^*$, and $c_j$ are rational numbers such that

$$\sum_{j=1}^k c_j f_j \otimes \cdots \otimes f_j \otimes (f_j \wedge (1-f_j)) = 0$$

in $Sym^{n-2}F_{\mathbb{Q}}^* \otimes \bigwedge^2 F_{\mathbb{Q}}^*$, then the function

$$z \mapsto \sum_{j=1}^k c_j P_n^{\text{mod}}(f_j(z))$$

is constant on $C$.

**Proof** This is done by Zagier in the proof of [Za, Proposition 3] for $C = \mathbb{P}^1_\mathbb{C}$, which works just as well for any curve as in the statement of the Lemma.

Because $B_0 = 0$, $B_1 = -\frac{1}{2}$ and $B_2 = \frac{1}{6}$, we have

$$P_{3,\text{Zag}}(z) = P_3^{\text{mod}}(z) - \frac{1}{6}\log^2|z|\log|1-z|. \tag{2.4}$$

Propositions 4.1, 5.1, Remark 5.2 and Theorem 5.3 of [dJ1] contain the following result.
**Theorem 2.3** Let $k$ be a number field and let $\sigma_1, \ldots, \sigma_r$ be all embeddings of $k$ into $\mathbb{C}$. Then the maps $\varphi^p_{(n)}$ and $\tilde{\varphi}^p_{(n)}$ exist without assumptions. They are injective for $p = 1$, and isomorhisms for $(p, n) = (1, 2)$ or $(1, 3)$. Moreover, the composition

$$H^1(\widetilde{\mathcal{M}}^\bullet_{(n)}(k)) \overset{\varphi^1_{(n)}}{\to} K^{(n)}_{2n-1}(k) \overset{\text{reg}}{\to} H^0_{\text{dR}}(\text{Spec}\,(k \otimes_\mathbb{Q} \mathbb{C}); \mathbb{R}(n-1))^+ = (\oplus_\sigma \mathbb{R}(n-1)_\sigma)^+.$$

is given by mapping $[x]_n$ to $\pm(n-1)! \oplus_\sigma P_n^{\text{mod}}(\sigma(x))_\sigma$.

Finally, Borel's theorem implies that for a number field $k$ the regulator on $K^{(n)}_{2n-1}(k)$ is an injection for $n \geq 2$.

We shall want the following theorem for the computation of the boundary map under localization.
**Theorem 2.4** We have a commutative diagram (up to sign)

$$\begin{array}{ccc}
K_n^{(n+1)}(X^n; \square^n) & \longrightarrow & K_n^{(n+1)}(X \times X_{\text{loc}}^{n-1}; \square^n) \\
\cong \downarrow & & \cong \downarrow \\
K_{n+1}^{(n+1)}(X^{n-1}; \square^{n-1}) & \longrightarrow & K_{n+1}^{(n+1)}(X_{\text{loc}}^{n-1}; \square^{n-1})
\end{array}$$

and the image of $\sum_j c_j g_j \otimes [f_j]_n$ in $K_n^{(n+1)}(X^n; \square^n)$ under the map

$$\varphi^2_{(n+1)} : H^2(\mathcal{M}^\bullet_{(n+1)}(F)) \to K_n^{(n+1)}(X^n; \square^n) \cong K_{2n}^{(n+1)}(F)$$

maps to $\pm \sum_j c_j[f_j]_n \cup g_j$ in $K_{n+1}^{(n+1)}(X_{\text{loc}}^{n-1}; \square^{n-1})$, modulo $(1+I)^* \tilde{\cup} K_n^{(n)}(X_{\text{loc}}^{n-2}; \square^{n-2})$.

**Proof** The proof is rather analogous to the proof of [dJ2, Proposition 3.2].

For computing the image under the regulator map we use integration. Because we shall be integrating forms on non–compact varieties, we need some results about the dependence of the result on the explicit representative chosen for a particular class. This problem was dealt with in Proposition 4.6 of [dJ2], which we now proceed to recall.

Let $Y$ be an algebraic variety of dimension $n$ with compactification $\overline{Y}$ such that the complement of $Y$ is given by $D$, a divisor with normal crossings. Suppose moreover that there is another divisor $D'$ on $\overline{Y}$ such that the union of $D$ and $D'$ is a divisor with normal crossings. We want to say something about the behaviour close to $D$ of forms on $Y$ that vanish on $Y \cap D'$. Suppose that locally in a compact subset of $\overline{Y}$, $D$ is given by $\prod_{i=1}^k x_i = 0$.

Let $r_i = |x_i|$, $\theta_i = \arg x_i$. We will consider differential forms $\beta$ on $Y$ that satisfy the following condition on the compact subset of the chosen neighbourhood of $D$.

$$\beta \text{ vanishes on } D' \cap Y \text{ and can be written as sums of products of } \log r_i, \, \mathrm{d}\log r_i, \, \mathrm{d}\theta_i, \text{ bounded functions on } Y \text{ and the restriction of } C^\infty\text{–forms on } \overline{Y} \text{ to } Y. \tag{2.5}$$

**Proposition 2.5** Let $Y$, $\overline{Y}$, $D$ and $D'$ be as above. Suppose that $\beta_1$ and $\beta_2$ are two closed $n$–forms on $Y$ as in (2.5) that represent the same class in relative de–Rham cohomology $H^n(Y; D'; \mathbb{R}(j))$. Let $\omega$ be a holomorphic or anti–holomorphic $n$–form on $\overline{Y}$, possibly with logarithmic poles along $D'$. Then

$$\int_Y \beta_1 \wedge \omega = \int_Y \beta_2 \wedge \omega.$$

We conclude this section with some remarks on orientations and standard integrals, to be used throughout the paper.

We shall always use the following orientations for the integrals involved: $t = x + iy$ the standard parameter on $\mathbb{A}^1 \subset \mathbb{P}^1$, then $\mathbb{P}^1$ or open parts have orientation given by $\mathrm{d}x \wedge \mathrm{d}y = \frac{-1}{2i}\mathrm{d}t \wedge \mathrm{d}\bar{t}$. On $(\mathbb{P}^1)^n$ or open parts, we use the product orientation given by $\frac{-1}{2i}\mathrm{d}t_1 \wedge \mathrm{d}\bar{t}_1 \wedge \ldots \wedge \frac{-1}{2i}\mathrm{d}t_n \wedge \mathrm{d}\bar{t}_n$. On $X_{S^1}^n = X^n \times S^1$ we take the product orientation of the above on $X^n$ with the standard counterclockwise orientation on $S^1$.

Using Stokes' theorem and the fact that for $c \in \mathbb{C}$, $\mathrm{d}\log \dfrac{t-c}{t-1} \wedge \mathrm{d}\log t = 0$, we find

$$\int_X \mathrm{d}i \arg \frac{t-c}{t-1} \wedge \mathrm{d}\,\log|t| = -\int_X \mathrm{d}\log\left|\frac{t-f}{t-1}\right| \wedge \mathrm{d}i \arg t = 2\pi i \log|c|$$

and

$$\int_X \mathrm{d}i \arg\left(\frac{t-f}{t-1}\right) \wedge \mathrm{d}i \arg t = -\int_X \mathrm{d}\log\left|\frac{t-f}{t-1}\right| \wedge \mathrm{d}\log|t| = 0,$$

hence

$$\int_X \mathrm{d}\log\left|\frac{t-c}{t-1}\right| \wedge \mathrm{d}\log \bar{t} = \int_X \mathrm{d}\log\left|\frac{t-c}{t-1}\right|^2 \wedge \mathrm{d}\log \bar{t} = 4\pi i \log|c|.$$

We also have the standard integral for $\rho_1$ a bump function around $t = 0$, i.e., $\rho_1 \equiv 1$ around $t = 0$ and $\rho_1 \equiv 0$ off $t = 0$,

$$\int_X \mathrm{d}(\rho_1(t)\mathrm{d}i \arg t) = -2\pi i.$$

Finally, we shall need the following integral. Let $h$ be a function on $\mathbb{P}^1$ with $h(\infty) - h(0) = 1$. Then

$$\int_X \mathrm{d}h(t) \wedge \mathrm{d}i \arg t = 2\pi i.$$

## 3 The regulator integral

Let $C$ be a smooth, proper, geometrically irreducible curve over the number field $k$, and let $g$ be its genus. Then $C_{\mathrm{an}}$, the associated complex manifold to $C \otimes_\mathbb{Q} \mathbb{C}$ is a disjoint union of $[k:\mathbb{Q}]$ complete algebraic curves $C_\tau$ over $\mathbb{C}$ of genus $g$, indexed by the embeddings of $k$ into $\mathbb{C}$. We fix an orientation on $C_{\mathrm{an}}$ such that the involution $\sigma$ given by complex conjugation on $\mathbb{C}$ in $C \otimes_\mathbb{Q} \mathbb{C}$ reverses the orientation. We also introduce the number $r$ defined by $r = [k:\mathbb{Q}] \cdot g$.

The goal of this section is to describe the regulator on the image of $\varphi^2_{(n+1)}$ inside $H^1_{\mathrm{dR}}(F; \mathbb{R}(n))^+$. We begin with some remarks on the cohomology groups of $C_{\mathrm{an}}$.

For $\tau : k \to \mathbb{C}$ an embedding, denote $C_\tau$ the curve obtained from $C$ by base change from $k$ to $C$ via $\tau$.

If $\tau$ is a real embedding, then $\sigma$ acts on $C_\tau$ by reversing its orientation. $H^1_{\mathrm{dR}}(C_\tau; \mathbb{C})$ is spanned (as a $\mathbb{C}$–vector space) by the holomorphic and the anti–holomorphic forms on $C_\tau$. Then the $\mathbb{R}$–vector space of holomorphic 1–forms $\omega$ on $C_\tau$ such that $\omega \circ \sigma = \overline{\omega}$ is given by $H^0(C_\mathbb{R}; \Omega) \cong \mathbb{R}^g$ where $C_\mathbb{R}$ is the base change from $k$ to $\mathbb{R}$ via $\tau$. On the other hand, by projecting

$$H^1_{\mathrm{dR}}(C_\tau; \mathbb{C}) \cong H^0(C_\mathbb{R}; \Omega) \oplus H^0(C_\mathbb{R}; \overline{\Omega}) \oplus iH^0(C_\mathbb{R}; \Omega) \oplus iH^0(C_\mathbb{R}; \overline{\Omega})$$

onto the $\mathbb{R}(n)$ and $\pm$ parts one checks easily that

$$H^1_{\mathrm{dR}}(C_\tau; \mathbb{R}(n))^+ = \pi_n H^0(C_\mathbb{R}; \Omega) \cong \mathbb{R}^g$$

because the forms remain independent after projection onto the real or imaginary parts.

The pairing here is given by mapping $(\pi_n \psi, \varphi)$ in $H^1_{\mathrm{dR}}(C_\tau; \mathbb{R}(n))^+ \times H^0(C_\mathbb{R}; \Omega)$ to

$$\int_{C_\tau} \pi_n \psi \wedge \overline{\varphi} = \frac{1}{2} \int_{C_\tau} \psi \wedge \overline{\varphi}$$

so that this pairing is non–degenerate is a consequence of the duality between holomorphic and anti–holomorphic forms on $C_\tau$.

If $\tau$ is not a real embedding, $\sigma$ acts on $C_\tau \coprod C_{\bar{\tau}}$. Then the holomorphic 1–forms $\omega$ such that $\omega \circ \sigma = \overline{\omega}$ are given by the pairs $(\omega, \overline{\omega \circ \sigma})$ with $\omega \in H^0(C_\tau; \Omega) \cong \mathbb{C}^g$. And $H^1_{\mathrm{dR}}(C_\tau \coprod C_{\bar{\tau}}; \mathbb{R}(n))^+$ is given by the pairs $(\psi, \overline{\psi \circ \sigma})$ with $\psi \in H^1_{\mathrm{dR}}(C_\tau; \mathbb{R}(n))$ which has $\mathbb{R}$–dimension $2g$. The pairing is given in this case by mapping $((\psi, \overline{\psi \circ \sigma}), (\omega, \overline{\omega \circ \sigma}))$ to $\int_{C_\tau} \psi \wedge \overline{\omega} + \int_{C_{\bar{\tau}}} \overline{\psi \circ \sigma} \wedge \omega \circ \sigma = 2 \int_{C_\tau} \psi \wedge \pi_{n+1} \omega$ because $\sigma$ reverses the orientation and $\overline{\psi} = (-1)^n \psi$. So this pairing is non–degerate because the pairing

$$H^1_{\mathrm{dR}}(C_\tau; \mathbb{R}(n)) \times H^1_{\mathrm{dR}}(C_\tau; \mathbb{R}(n+1)) \to H^2_{\mathrm{dR}}(C_\tau; \mathbb{R}(1)) \cong \mathbb{R}(1)$$

is non–degenerate and the projection $H^0(C_\tau; \Omega) \to H^1_{\mathrm{dR}}(C_\tau; \mathbb{R}(n+1))$ given by mapping $\omega$ to $\pi_{n+1} \omega$ is an isomorphism.

This shows that $H^1_{\mathrm{dR}}(C; \mathbb{R}(n))^+ \cong \mathbb{R}^r$ with $r = [k:\mathbb{Q}] \cdot g$ as above. We also have the following duality.

**Remark 3.1** The holomorphic forms $\omega$ such that $\omega \circ \sigma = \overline{\omega}$ form the dual of $H^1_{\mathrm{dR}}(C; \mathbb{R}(n))^+$ under the pairing

$$H^1_{\mathrm{dR}}(C; \mathbb{R}(n))^+ \times <\omega> \to \mathbb{R}(1)$$

defined by sending $(\psi, \omega)$ to $\int_C \psi \wedge \overline{\omega}$.

We use this duality in Proposition 3.2 below in order to give the Beilinson regulator in terms of integrals.

Write $\pm$ or $\mp$ and in formulas read either the top or the bottom in all places. The involution $\sigma$ acts also on $H_1(C;\mathbb{Q})$, so this space splits into a $+$– and a $-$–part as well. From the pairing with $H^1(C;\mathbb{Q})$ one deduces that both pieces have $\mathbb{Q}$–dimension $r$, as $H_1(\mathbb{Q})^\pm$ is perpendicular to $H^1(\mathbb{Q})^\mp$. Let $\{s_{1,\pm},\ldots,s_{r,\pm}\}$ be a basis of $H_1(C;\mathbb{Q})^\pm$, and let $\{s^*_{1,\pm},\ldots,s^*_{r,\pm}\}$ in $H^1(C;\mathbb{Q})$ be its dual base, so that $\int_{s_{m,\pm}} s^*_{k,\pm} = \delta_{mk}$. Let $T^\pm_{k,l} = \left(\int_{s_{k,\pm}} \omega_l\right)$, so $\omega_l = \sum_k \left(T^+_{k,l} s^*_{k,+} + T^-_{k,l} s^*_{k,-}\right)$. Write $\mathbb{R}(+)$ for $\mathbb{R}(0)$ and $\mathbb{R}(-)$ for $\mathbb{R}(1)$, and similarly for $\pi_\pm$. As $\overline{\omega_l} = \omega_l \circ \sigma$, we get $\int_{s_{k,\pm}} \overline{\omega_l} = \int_{s_{k,\pm}} \omega_l \circ \sigma = \int_{\sigma(s_{k,\pm})} \omega_l = \pm \int_{s_{k,\pm}} \omega_l$, and hence $T^\pm$ has entries in $\mathbb{R}(\pm)$. Therefore $T^\pm_{kl} = \int_{s_{k,\pm}} \pi_\pm \omega_l$. In particular, $\pi_+ \omega_l = \sum_k T^+_{k,l} s^*_{k,+}$ and $\pi_- \omega_l = \sum_k T^-_{k,l} s^*_{k,-}$.

**Proposition 3.2** Suppose the elements $\alpha_1,\ldots,\alpha_r$ of $K^{(n+1)}_{2n}(C)$ under the Beilinson regulator get mapped to $\psi_1,\ldots,\psi_r$ in $H^1_{\mathrm{dR}}(C;\mathbb{R}(n))^+$. Let $\omega_1,\ldots,\omega_r$ and $T^\pm$ be as above, and let $R_{k,l} = \int_C \psi_k \wedge \overline{\omega_l}$. Then the Beilinson regulator of $\alpha_1,\ldots,\alpha_r$ is given by

$$c_{n+1} = \frac{\det(R)}{(2\pi i)^{nr} \det(T^\pm)}.$$

where we take $-$ for $n$ even, and $+$ for $n$ odd.

**Proof** In this proof, let us write $\pm$ and $\mp$ where we mean that we take the top for $n$ even, the bottom for $n$ odd. Note that $\psi_k \circ \sigma = \pm \psi_k$ if $\psi_k \in H^1_{\mathrm{dR}}(C;\mathbb{R}(n))^+$. Therefore we can define the $\mathbb{R}(n)$–valued matrix $M$ by $\psi_k = \sum_m M_{k,m} s^*_{m,\pm}$. Then by definition, $c_{n+1} = (2\pi i)^{-nr} \det(M)$. Note that as $\sigma$ reverses the orientation, $\omega_l \circ \sigma = \overline{\omega_l}$ and $\psi_k \circ \sigma = \overline{\psi_k} = (-1)^n \psi_k$,

$$R_{k,l} = \int_C \psi_k \wedge \overline{\omega_l} = -\int_C \psi_k \circ \sigma \wedge \overline{\omega_l \circ \sigma} = -\int_C \overline{\psi_k} \wedge \omega_l = -\overline{R_{k,l}}.$$

Therefore $R_{k,l}$ is purely imaginary, and

$$\begin{aligned} R_{k,l} &= -\int_C \psi_k \wedge \pi_\mp \omega_l \\ &= -\sum_n T^\mp_{n,l} \int_C \psi_k \wedge s^*_{n,\mp} \\ &= -\sum_{m,j} M_{j,m} T^\mp_{j,l} \int_C s^*_{m,\pm} \wedge s^*_{n,\mp} \\ &= -\sum_{m,j} M_{j,m} A_{m,j} T^\mp_{j,l} \end{aligned}$$

with $A_{m,j} = \int_C s^*_{m,\pm} \wedge s^*_{n,\mp}$. As $\det(A)$ expresses the non–degeneracy of the pairing $H^1(C;\mathbb{Q}) \times H^1(C;\mathbb{Q}) \to H^2(C;\mathbb{Q})$, it is an element of $\mathbb{Q}^*$. Hence taking determinants we find that (up to a factor in $\mathbb{Q}^*$) $\det(R) = \det(M) \det(T^\mp)$. So we get that the regulator $c_{n+1}$ of $\alpha_1,\ldots,\alpha_r$ is given by $c_{n+1} = \dfrac{\det(R)}{(2\pi i)^{nr} \det(T^\mp)}$.

**Notation** If $\omega$ is a holomorphic 1–form on $C$ such that $\omega \circ \sigma = \omega$, we call the map

$$H^1_{\mathrm{dR}}(C;\mathbb{R}(n))^+ \to \mathbb{R}(1)$$

given by
$$\psi \mapsto \int_C \psi \wedge \overline{\omega}$$
the regulator integral associated to $\omega$. We shall use the same terminology if we precede this with the regulator map $K_{2n}^{(n+1)}(C) \to H_{\mathrm{dR}}^1(C; \mathbb{R}(n))^+$.

The regulator integral has the advantage that it can be factored over larger groups than just $K_{2n}^{(n+1)}(C)$:

**Proposition 3.3** Let $\omega$ be a holomorphic 1–form on $C$. Then the regulator integral $K_{2n}^{(n+1)}(C) \to H_{\mathrm{dR}}^1(C; \mathbb{R}(n))^+ \to \mathbb{R}(1)$ associated to $\omega$ extends naturally over the maps
$$K_{2n}^{(n+1)}(C) \to K_{2n}^{(n+1)}(F) \to K_{2n}^{(n+1)}(F)/K_{2n-1}^{(n)}(k) \cup F_{\mathbb{Q}}^*.$$

**Proof** Let $\beta$ be a class in $H_{\mathrm{dR}}^1(F; \mathbb{R}(n))^+ = \lim_{U \subset C} H_{\mathrm{dR}}^1(U; \mathbb{R}(n))^+$. Then using the fact that $H_{\mathrm{dR}}^1(U; \mathbb{C})$ can be computed using forms with logarithmic singularities, one sees that $\beta$ has a representative $\psi$ as in (2.5), and we extend the map by mapping $\beta$ to $\int_C \psi \wedge \overline{\omega}$. Proposition 2.5 shows that this does not depend on our choice of $\psi$, hence is well defined.

As for the last map, note that the regulator of $\alpha \cup f$ for $\alpha \in K_{2n-1}^{(n)}(k)$ and $f \in F$ is given by $(0, c) \cup (\mathrm{d}\log f, \log|f|) = (0, c\mathrm{d}i \arg f)$ in Deligne cohomology, hence by $d\mathrm{d}i \arg f$ in $H_{\mathrm{dR}}^1$. Then the regulator integral becomes $\int_C c\mathrm{d}i \arg f \wedge \overline{\omega} = \int_C c\mathrm{d}\log|f| \wedge \overline{\omega} = 0$ as one easily checks using Stokes' theorem.

The rest of the section is devoted to rewriting the integrals occurring in Proposition 3.2 on the image of $H^2(\mathcal{M}_{(n+1)}^\bullet(F))$. We shall in fact prove the following theorem.

**Theorem 3.4** Let $\sum_j c_j g_j \otimes [f_j]_n$ be an element of $H^2(\mathcal{M}_{(n+1)}^\bullet(F))$. Then the regulator integral
$$H^2(\mathcal{M}_{(n+1)}^\bullet(F)) \longrightarrow K_{2n}^{(n+1)}(F) \xrightarrow{\mathrm{reg}} H_{\mathrm{dR}}^1(F; \mathbb{R}(n)) \xrightarrow{\int_C \cdots \wedge \overline{\omega}} \mathbb{R}(1)$$
is given by mapping $g \otimes [f]_n$ to
$$\pm 2^n \int_C \log|g| \log^{n-1}|f| \mathrm{d}\log|1 - f| \wedge \overline{\omega}$$

One can replace this with
$$\pm 2^n \frac{n}{n+1} \int_C \log|g| \log^{n-2}|f| \left(\log|1 - f|\mathrm{d}\log|f| - \log|f|\mathrm{d}\log|1 - f|\right) \wedge \overline{\omega}, \qquad (3.1)$$
which factors through the map $H^2(\mathcal{M}_{(n+1)}^\bullet(F)) \to H^2(\widetilde{\mathcal{M}}_{(n+1)}^\bullet(F))$, i.e., vanishes on symbols $<[f]_n + (-1)^n [1/f]_n> \otimes g$.

**Remark 3.5** If the curve is an elliptic curve $E$, the integrals occurring in Theorem 3.4 can be rewritten using Fourier transformation. This gives expressions for the regulator integral in terms of non–classical Eisenstein series (see e.g., [G3, Theorem 3.4] for the case $n = 2$ and [G1, Theorem 5.8] for arbitrary $n$). It seems that for $n = 2$ those Eisenstein series were first considered by Deninger in [De].

The rest of this section is devoted to the proof of Theorem 3.4. We begin with making the isomorphism $H^{n+1}_{\mathrm{dR}}(X^n_C; \square^n; \mathbb{R}(n))^+ \cong H^1_{\mathrm{dR}}(C; \mathbb{R}(n))^+$ explicit. Namely, let $h$ be a real–valued function on $\mathbb{P}^1$ such that $h(\infty) - h(0) = 1$. Then the isomorphism

$$H^{n-1}_{\mathrm{dR}}(Y; \mathbb{R}(m)) \tilde{\to} H^n_{\mathrm{dR}}(X_Y; \square; \mathbb{R}(m))$$

is given by sending $\psi$ to $\psi \wedge \mathrm{d}h$, and similarly for the $+$–parts if $h$ is symmetric with respect to complex conjugation on $\mathbb{P}^1_\mathbb{C}$. As

$$\int_X \mathrm{d}h \wedge \mathrm{d}\log \bar{t} = -\int_X \mathrm{d}h \wedge \mathrm{d}i\arg t = -2\pi i$$

this means that for $\psi$ in $H^1_{\mathrm{dR}}(C; \mathbb{C})$

$$\pm \int_C \psi \wedge \bar{\omega}$$
$$= (2\pi i)^{-n} \int_{X^n_C} \psi \wedge \mathrm{d}h(t_1) \wedge \ldots \wedge \mathrm{d}h(t_n) \wedge \mathrm{d}\log\bar{t_1} \wedge \ldots \wedge \mathrm{d}\log\bar{t_n} \wedge \bar{\omega}$$
$$= (-2\pi i)^{-n} \int_{X^n_C} \psi \wedge \mathrm{d}h(t_1) \wedge \ldots \wedge \mathrm{d}h(t_n) \wedge \mathrm{d}i\arg t_1 \wedge \ldots \wedge \mathrm{d}i\arg t_n \wedge \bar{\omega}.$$

Now let $\alpha = \sum_j c_j [f_j]_n \otimes g_j$ be an element of $H^2(\mathcal{M}^\bullet_{(n+1)}(F))$, with image $\psi$ in $H^1_{\mathrm{dR}}(F; \mathbb{R}(n))^+$ under the composition

$$H^2(\mathcal{M}^\bullet_{(n+1)}(F)) \to K^{(n+1)}_{2n}(F) \to H^1_{\mathrm{dR}}(F; \mathbb{R}(n))^+.$$

From Theorem 2.4 we have a commutative diagram

$$\begin{array}{ccc}
K^{(n+1)}_n(X^n; \square^n) & \longrightarrow & K^{(n+1)}_n(X \times X^{n-1}_{\mathrm{loc}}; \square^n) \\
\downarrow \cong & & \downarrow \cong \\
K^{(n+1)}_{n+1}(X^{n-1}; \square^{n-1}) & \longrightarrow & K^{(n+1)}_{n+1}(X^{n-1}_{\mathrm{loc}}; \square^{n-1})
\end{array}$$

and the image of $\alpha = \sum_j c_j[f_j]_n \otimes g_j$ under $\varphi^2_{(n+1)}$ in $K^{(n+1)}_n(X^n; \square^n)$ maps to $\alpha' = \pm \sum_j c_j[f_j]_n \cup g_j$ in $K^{(n+1)}_{n+1}(X^{n-1}_{\mathrm{loc}}; \square^{n-1})$, modulo $(1+I)^* \tilde{\cup} K^{(n)}_n(X^{n-2}_{\mathrm{loc}}; \square^{n-2})$. One has the corresponding diagram in Deligne cohomology, which is equal to de deRham cohomology in all cases, so that we have $\mathrm{reg}(\alpha)$ in $H^{n+1}_{\mathrm{dR}}(X^n; \square^n; \mathbb{R}(n))^+$, which corresponds to a form $\psi'$ in $H^n_{\mathrm{dR}}(X^{n-1}; \square^{n-1}; \mathbb{R}(n))^+$, which in turn maps to $\mathrm{reg}(\alpha')$ in $H^n_{\mathrm{dR}}(X^{n-1}_{\mathrm{loc}}; \square^{n-1}; \mathbb{R}(n))^+$. Therefore, where $\bar{\omega}$ is any holomorphic 1–form on $C$, and $\Theta = \mathrm{d}\log\bar{t_1} \wedge \ldots \wedge \mathrm{d}\log\bar{t_{n-1}}$, we want to compute

$$\pm \int_C \psi \wedge \bar{\omega} = \frac{1}{(2\pi i)^{n-1}} \int_{X^{n-1} \times C} \psi' \wedge \Theta \wedge \bar{\omega} = \frac{1}{(2\pi i)^{n-1}} \int_Z \mathrm{reg}(\alpha') \wedge \Theta \wedge \bar{\omega}$$

where $Z$ is a localization of $X^{n-1}_C$ and $\mathrm{reg}(\alpha')$ is a form satisfying (2.5) on a suitable compatification of $Z$, $\bar{Z}$. We can obtain $\bar{Z}$ from $(\mathbb{P}^1)^{n-1} \times C$ by repeatingly blowing up, obtaining a $Z$ which has a Zariski open part isomorphic to $X_{U,\text{loc}}^{n-1}$ for some Zariski open $U$ of $C$. To see that the last equality sign is true, we have to check that the conditions of Proposition 2.5 apply, i.e., that on $\overline{Z}$ reg($\alpha'$) can be written as in (2.5) and that (the pullback to $\overline{Z}$ of) $\operatorname{d} \log \overline{t_1} \wedge \ldots \wedge \operatorname{d} \log \overline{t_{n-1}} \wedge \overline{\omega}$ has poles of order one along the strict transform of $\square^{n-1}$ and no poles elsewhere. This one easily checks explicitly, cf. the computations on page 608 of [G–H].

Because we can take $\overline{Z}$ to be a blowup of $(\mathbb{P}^1)_C^{n-1}$, which is isomorphic to $(\mathbb{P}^1)_U^{n-1}$ for some Zariski open part $U$ of $C$, we can compute simply on $X_{U,\text{loc}}^{n-1}$ without changing the value of the integral

As the regulator of $[f]_n \cup g$ is given by the product $(\omega_n, \varepsilon_n) \cup (\operatorname{d}\log g, \log |g|) = (\omega_n \wedge \operatorname{d}\log g, \varepsilon_n \wedge \operatorname{d} i \arg g + (-1)^n \log |g| \pi_n \omega_n)$, we find that the regulator integral is given by

$$(2\pi i)^{1-n} \sum_j c_j \int_{\overline{Z}} (\varepsilon_n(f_j) \wedge \operatorname{d} i \arg g_j + (-1)^n \log |g_j| \pi_n \omega_n(f_j)) \wedge \Theta \wedge \overline{\omega}$$

$$=(2\pi i)^{1-n} \sum_j c_j \int_{\overline{Z}} (\varepsilon_n(f_j) \wedge \operatorname{d} \log |g_j| + (-1)^n \log |g_j| \pi_n \omega_n(f_j)) \wedge \Theta \wedge \overline{\omega}$$

as $\operatorname{d}(i \arg g - \log |g|) \wedge \overline{\omega}(f) = 0$ on $X^{n-1} \times C$. Adding

$$0 = (-1)^n (2\pi i)^{1-n} \sum_j c_j \int_{\overline{Z}} \operatorname{d}(\log |g| \varepsilon_n(f_j)) \wedge \Theta \wedge \overline{\omega}$$

we obtain

$$(-1)^n (2\pi i)^{1-n} \sum_j c_j \int_{\overline{Z}} \log |g_j| \, \omega_n(f_j) \wedge \Theta \wedge \overline{\omega}$$

Remembering that

$$\omega_n(f) = (-1)^{n-1} \operatorname{d} \log \frac{t_1 - f}{t_1 - 1} \wedge \ldots \wedge \operatorname{d} \log \frac{t_{n-1} - f}{t_{n-1} - 1} \wedge \operatorname{d} \log(1 - f)$$

we can compute this as

$$\pm (2\pi i)^{1-n} \sum_j c_j \int_C \left( \int_X \operatorname{d} \log \frac{t - f_j}{t - 1} \wedge \operatorname{d} \log \overline{t} \right)^{n-1} \log |g_j| \operatorname{d} \log(1 - f_j) \wedge \overline{\omega}$$

$$= \pm \sum_j c_j 2^{n-1} \int_C \log |g_j| \log^{n-1} |f_j| \operatorname{d} \log(1 - f_j) \wedge \overline{\omega}$$

$$= \pm \sum_j c_j 2^n \int_C \log |g_j| \log^{n-1} |f_j| \operatorname{d} \log |1 - f_j| \wedge \overline{\omega}$$

because $\int_X \log \frac{t-f}{t-1} \wedge \operatorname{d} \log \overline{t} = 2\pi i \log^2 |f|$ and $\operatorname{d} \log \overline{(1-f)} \wedge \overline{\omega} = 0$ on $C$.

In order to rewrite this in terms of

$$\int_C \log |g| \log^{n-2} |f| \left( \log |1 - f| \operatorname{d} \log |f| - \log |f| \operatorname{d} \log |1 - f| \right) \wedge \overline{\omega}$$

we add
$$-\sum_j c_j \frac{1}{n+1} \int_C \mathrm{d}\left(\log|g_j|\log^{n-1}|f_j|\log|1-f_j|\right) \wedge \overline{\omega}$$
and
$$\sum_j c_j \frac{1}{n+1} \int_C \log^{n-2}|f_j|\log|1-f_j|\left(\log|f_j|\mathrm{d}\log|g_j| - \log|g_j|\mathrm{d}\log|f_j|\right) \wedge \overline{\omega}.$$

Note that this does not change the value of the integral, as the first term vanishes by Stokes' theorem, and the second because we take sums of terms corresponding to an element in $H^2(\mathcal{M}^\bullet_{(n+1)}(F))$, and the form vanishes identically after summing up the terms. This yields the integral given in (3.1). One checks immediately by writing it out that this vanishes on $([f]_n + (-1)^n[1/f]_n) \otimes g$.

**Remark 3.6** For $n+1 = 4$, i.e., the case $K_6^{(4)}(F)$, we get those results without any assumptions. The map $\varphi^2_{(4)}$ from $H^2(\mathcal{M}^\bullet_{(4)}(F))$ to $K_6^{(4)}(F)/N$ (with $N$ generated by $K_4^{(2)}(F) \cup K_2^{(2)}(F)$, see Section 2) exists without assumptions, and the regulator factors through this as the regulator map to Deligne cohomology respects the product structure, and vanishes on $K_4^{(2)}(F)$. Hence it factors through this quotient to give us
$$H^2(\mathcal{M}^\bullet_{(4)}(F)) \to H^2_\mathcal{D}(F, \mathbb{R}(3))^+ \cong H^1_{\mathrm{dR}}(F, \mathbb{R}(2))^+$$
as in Theorem 3.4.

## 4 The boundary under localization

Because $K_{2n}(k(x))$ is torsion as $k(x)$ is a number field, the localization sequence takes the form
$$0 \to K^{(n+1)}_{2n}(C) \to K^{(n+1)}_{2n}(F) \to \coprod_{x \in C^{(1)}} K^{(n)}_{2n-1}(k(x)) \to \dots.$$

This section is devoted to the computation of the boundary map on the image in $K^{(n+1)}_{2n}(F)$ of $H^2(\widetilde{\mathcal{M}}^\bullet_{(n)}(F))$ or $H^2(\mathcal{M}^\bullet_{(n)}(F))$ in this localization sequence for $n = 3$. The method chosen probably works for all $n \geq 2$ (with the case $n = 2$ already done in [dJ2]), but at some stage the combinatorics become too complicated in general and we restrict ourselves to the case $n = 3$.

Recall that in [dJ2, Corollary 5.4] it was proved that we have a commutative diagram (up to sign and up to $\partial\left(K_3^{(2)}(k) \cup F_\mathbb{Q}^*\right)$ in the lower right hand corner)

$$\begin{CD}
H^2(\mathcal{M}^\bullet_{(3)}(F)) @>{\varphi^2_{(3)}}>> K_4^{(3)}(F) \\
@V{2\delta}VV @VV{\partial}V \\
\coprod_{x \in C^{(1)}} H^1(\widetilde{\mathcal{M}}^\bullet_{(2)}(k(x))) @>{\sim}>> \coprod_{x \in C^{(1)}} K_3^{(2)}(k(x))
\end{CD}$$

Note that the lower horizontal arrow is an isomorphism by Theorem 2.3. For generalizing this to $n = 3$, we need some preliminary results. The following result was proved in [dJ2] for $n = 2$.

**Proposition 4.1** There is a map

$$\delta : \quad \widetilde{\mathcal{M}}^\bullet_{(n)}(F) \to \coprod_{x \in C^{(1)}} \widetilde{\mathcal{M}}^\bullet_{(n-1)}(k(x))[1]$$

given by

$$\delta_x([f]_{n-k} \otimes g_1 \wedge \ldots \wedge g_k) = \mathrm{sp}_{n-k,x}([f]_{n-k}) \otimes \partial_{k,x}(g_1 \wedge \ldots \wedge g_k)$$

for $k = 1, \ldots, n-1$, and

$$\delta_x(g_1 \wedge \ldots \wedge g_n) = -\partial_{n,x}(g_1 \wedge \ldots \wedge g_n)$$

for $k = n$, where $\mathrm{sp}_{n-k,x}([f]_{n-k}) = [f(x)]_{n-k}$ if $f(x) \neq 0$ or $\infty$, 0 otherwise, and $\partial_{k,x}$ the unique map from $\bigwedge^k F^*_\mathbb{Q}$ to $\bigwedge^{k-1} k(x)^*_\mathbb{Q}$ determined by

$$\pi_x \wedge u_1 \wedge \ldots \wedge u_{k-1} \mapsto u_1(x) \wedge \cdots \wedge u_{k-1}(x)$$
$$u_1 \wedge \ldots \wedge u_k \mapsto 0$$

if all $u_i$ are units at $x$ and $\pi_x$ is a uniformizer at $x$. This gives rise to maps

$$\delta : \quad H^m(\widetilde{\mathcal{M}}^\bullet_{(n)}(F)) \to \coprod_{x \in C^{(1)}} H^{m-1}(\widetilde{\mathcal{M}}^\bullet_{(n-1)}(k(x)))$$

**Remark 4.2** We get induced maps

$$\delta : \quad \mathcal{M}^\bullet_{(n)}(F) \to \coprod_{x \in C^{(1)}} \widetilde{\mathcal{M}}^\bullet_{(n-1)}(k(x))[1]$$

and

$$\delta : \quad H^m(\mathcal{M}^\bullet_{(n)}(F)) \to \coprod_{x \in C^{(1)}} H^{m-1}(\widetilde{\mathcal{M}}^\bullet_{(n-1)}(k(x)))$$

by composing the natural projection $\mathcal{M}^\bullet_{(n)}(F) \to \widetilde{\mathcal{M}}^\bullet_{(n)}(F)$ with $\delta$.

**Proof** Let $x \in C$ be a closed point in our curve over the number field $k$. Fix a uniformizer $\pi_x$ around $x$. We shall in fact construct the map $\mathrm{sp}_{n,x} : M_{(n)}(F) \to \widetilde{M}_{(n)}(k(x))$, and then observe that it factors through the projection $M_{(n)}(F) \to \widetilde{M}_{(n)}(F)$.

Assume we have a map $\mathrm{sp}_{n-1,x} : M_{(n-1)}(F) \to \widetilde{M}_{(n-1)}(k(x))$ given by mapping $[f]_{n-1}$ to $[f(x)]_{n-1}$ if $f(x) \neq 0$ or $\infty$, and 0 otherwise. This was done for $n - 1 = 2$ in the proof of [dJ2, Prop. 5.1], and is the case where one has to work with $\widetilde{M}_{(n-1)}(k(x))$ rather than $M_{(n-1)}(k(x))$. We then have a diagram

$$\begin{array}{ccc} \mathbb{Q}[F^*] & \longrightarrow & M_{(n-1)}(F) \otimes F^*_\mathbb{Q} \\ \downarrow & & \downarrow \\ \widetilde{M}_{(n)}(k(x)) & \longrightarrow & \widetilde{M}_{(n-1)}(k(x)) \otimes k^*_\mathbb{Q} \end{array}$$

where $\mathbb{Q}[F^*]$ is the free $\mathbb{Q}$–vector space on elements of $F^*$, the vertical maps are $f \mapsto \mathrm{sp}_{n,x}([f(x)]_n)$ and $[f]_n \otimes g \mapsto \mathrm{sp}_{n-1,x}([f]_{n-1}) \otimes \overline{g(x)}$, with $\overline{g(x)} = g\pi_x^{-\mathrm{ord}_x(g)}|_x$, $\pi_x$ a uniformizer at $x$. It is obvious that the diagram commutes. To show that it factors through $M_{(n)}(F)$ observe that if $\alpha$ goes to zero in $M_{(n)}(F)$ then $\mathrm{sp}_{n,x}(\alpha)$ defines an element in $K^{(n)}_{2n-1}(k(x))$. As $k(x)$ is a number field, we can verify that it is zero by computing the regulator map, given by Theorem 2.3. To this end, consider all embeddings of $k(x)$ into $\mathbb{C}$, i.e., tensor the curve $C$ over $\mathbb{Q}$ with $\mathbb{C}$. Then we have that $P_n^{\mathrm{mod}}(\alpha)$ is constant, see Lemma 2.2. Specialize to a point $y$ where it can be done directly (which means that $y$ should lie in some Zariski open part, see Section 2), we find 0, so the regulator vanishes. Then we use that the regulator does not change if we replace $y$ with $x$, and continuity, to see that $P_n^{\mathrm{mod}}(\mathrm{sp}_{n,x}(\alpha)) = P_n^{\mathrm{mod}}(\mathrm{sp}_{n,y}(\alpha)) = 0$ because $P_n^{\mathrm{mod}}$ is continous at 0 and $\infty$ and has value 0. Hence $\mathrm{sp}_{n,x}(\alpha) = 0$ in $\widetilde{M}_{(n)}(k(x))$.

This map gives us the map $\mathrm{sp}_{n,x} : M_{(n)}(F) \to \widetilde{M}_{(n)}(k(x))$, obviously factoring through $\widetilde{M}_{(n)}(F)$. It is then easy to check that we get maps of complexes

$$\mathcal{M}^\bullet_{(n)}(F) \longrightarrow \widetilde{\mathcal{M}}^\bullet_{(n)}(F) \xrightarrow{\delta_x} \widetilde{\mathcal{M}}^\bullet_{(n-1)}(k(x))[1]$$

with $\delta_x$ given by mapping $[f]_{n-k} \otimes g_1 \wedge \ldots \wedge g_k$ to $\mathrm{sp}_{n-k,x}([f(x)]_{n-k}) \otimes \partial_{k,x}$ for $k = 1, \ldots, n-1$ and $-\partial_{n,x}(g_1 \wedge \ldots \wedge g_n)$ for $k = n$, with $\partial_{k,x}$ the unique map from $\bigwedge^k F^*_\mathbb{Q}$ to $\bigwedge^{k-1} k(x)^*_\mathbb{Q}$ given by

$$\pi_x \wedge u_1 \wedge \ldots \wedge u_{k-1} \mapsto u_1(x) \wedge \cdots \wedge u_{k-1}(x)$$
$$u_1 \wedge \ldots \wedge u_k \mapsto 0$$

if all $u_i$ are units at $x$. From this we get maps

$$\delta_x : H^m(\mathcal{M}^\bullet_{(n)}(F)) \to H^{m-1}(\widetilde{\mathcal{M}}^\bullet_{(n-1)}(k(x)))$$

as claimed in the proposition.

We can now introduce the complexes $\mathcal{M}^\bullet_{(n+1)}(C)$ and $\widetilde{\mathcal{M}}^\bullet_{(n+1)}(C)$, by defining them to be the total complexes of the double complexes

$$\begin{array}{ccccccccc}
0 & \to & M_{(n+1)}(F) & \xrightarrow{d} & M_{(n)}(F) \otimes F^*_\mathbb{Q} & \xrightarrow{d} & M_{(n-1)}(F) \otimes \bigwedge^2 F^*_\mathbb{Q} & \xrightarrow{d} & \cdots \\
 & & \downarrow & & \delta \downarrow & & \delta \downarrow & & \\
 & & 0 & \to & \coprod \widetilde{M}_{(n)}(k(x)) & \xrightarrow{d} & \coprod \widetilde{M}_{(n-2)}(k(x)) \otimes k(x)^*_\mathbb{Q} & \xrightarrow{d} & \cdots
\end{array}$$

and

$$\begin{array}{ccccccccc}
0 & \to & \widetilde{M}_{(n+1)}(F) & \xrightarrow{d} & \widetilde{M}_{(n)}(F) \otimes F^*_\mathbb{Q} & \xrightarrow{d} & \widetilde{M}_{(n-1)}(F) \otimes \bigwedge^2 F^*_\mathbb{Q} & \xrightarrow{d} & \cdots \\
 & & \downarrow & & \delta \downarrow & & \delta \downarrow & & \\
 & & 0 & \to & \coprod \widetilde{M}_{(n)}(k(x)) & \xrightarrow{d} & \coprod \widetilde{M}_{(n-2)}(k(x)) \otimes k(x)^*_\mathbb{Q} & \xrightarrow{d} & \cdots
\end{array}$$

where both coboundaries have degree 1 and the total complexes are cohomological complexes with $\mathcal{M}^\bullet_{(n+1)}(F)$ and $\widetilde{\mathcal{M}}^\bullet_{(n+1)}(F)$ in degree 1. There is an obvious inclusion of $H^2(\mathcal{M}^\bullet_{(n+1)}(C))$ into $H^2(\mathcal{M}^\bullet_{(n+1)}(F))$, and similarly we have an inclusion of $H^2(\widetilde{\mathcal{M}}^\bullet_{(n+1)}(C))$

into $H^2(\widetilde{\mathcal{M}}^\bullet_{(n+1)}(F))$, so that the maps $\varphi^2_{(n+1)}$ resp. $\tilde{\varphi}^2_{(n+1)}$ obviously extend to maps on the cohomology of those complexes.

**Corollary 4.3** The map $\varphi^2_{(3)}$ maps $H^2(\widetilde{\mathcal{M}}^\bullet_{(3)}(C))$ to $K^{(3)}_4(C) + K^{(2)}_3(k) \cup F^*_\mathbb{Q}$ inside $K^{(3)}_4(F)$.

**Remark 4.4** The exact sequence

$$0 \to K^{(3)}_4(C) \to K^{(3)}_4(C) + K^{(2)}_3(k) \cup F^*_\mathbb{Q} \to \partial\left(K^{(2)}_3(k) \cup F^*_\mathbb{Q}\right)$$

is split. Namely, because $K^{(2)}_3(k) \cup k = 0$ as $K^{(3)}_4(k) = 0$ (remember that $k$ is a number field), the cup product factors through $K^{(2)}_3(k) \otimes F^*_\mathbb{Q}/k^*_\mathbb{Q}$. The boundary map factors through this as well, and is given by $\partial \alpha \cup f = -\partial f \cup \alpha = -\text{div}(f) \cup \alpha$. Because the divisor map is injective on $F^*_\mathbb{Q}/k^*_\mathbb{Q}$, the boundary map is injective as well. The corresponding result holds for $K^{(n+1)}_{2n}(C)$ and $K^{(n)}_{2n-1}(k) \cup F^*_\mathbb{Q}$.

**Corollary 4.5** The map $\varphi^2_{(3)} : H^2(\mathcal{M}^\bullet_{(3)}(C)) \to K^{(3)}_4(C) + K^{(2)}_3(k) \cup F^*_\mathbb{Q}$ can be lifted to a map $\varphi^2_{(3)} : H^2(\mathcal{M}^\bullet_{(3)}(C)) \to K^{(3)}_4(C)$ by changing $\varphi^2_{(3)}(\alpha)$ with elements in $K^{(2)}_3(k) \cup F^*_\mathbb{Q}$.

We serve the main course in this section:

**Theorem 4.6** For $n = 3$, the diagram

$$\begin{array}{ccc}
H^2(\mathcal{M}^\bullet_{(4)}(F)) & \xrightarrow{\varphi^2_{(4)}} & K^{(4)}_6(F)/K^{(2)}_4(F) \cup K^{(2)}_2(F) \\
\downarrow{3\delta} & & \downarrow{\partial} \\
\coprod_{x \in C^{(1)}} H^1(\widetilde{\mathcal{M}}^\bullet_{(3)}(k(x))) & \xrightarrow{\sim} & \coprod_{x \in C^{(1)}} K^{(3)}_5(k(x))
\end{array}$$

commutes up to sign and up to $\partial\left(K^{(3)}_5(k) \cup F^*_\mathbb{Q}\right)$ in the lower right hand corner. (Note that the lower isomorphism is part of Theorem 2.3.)

**Remark 4.7** Note that from the localization sequence

$$0 \to K^{(2)}_4(C) \to K^{(2)}_4(F) \to K^{(1)}_3(k(x))$$

we get an isomorphism $K^{(2)}_4(F) \cong K^{(2)}_4(C)$ as $K^{(1)}_3(k(x)) = 0$. Hence on $K^{(2)}_4(F) \cup K^{(2)}_2(F)$ $\partial_x$ is given by mapping $\alpha \cup \beta$ to the cup product of $\alpha(x)$ with the boundary at $x$ of $\beta$, i.e., zero as $\alpha(x) \in K^{(2)}_4(k(x)) = 0$ as $k(x)$ is a number field. (Of course that would be zero for any field for which the Beilinson–Soulé conjecture holds.) So in fact $K^{(2)}_4(F) \cup K^{(2)}_2(F) \subset K^{(4)}_6(C)$.

**Corollary 4.8** The map $\varphi^2_{(4)}$ maps $H^2(\widetilde{M}_{(4)}(C))$ to $K^{(4)}_6((C))/K^{(4)}_2(F) \cup K^{(2)}_2(F) + K^{(3)}_5(k) \cup F^*_\mathbb{Q}$ inside $K^{(4)}_6(F)/K^{(2)}_4(F) \cup K^{(2)}_2(F)$.

**Corollary 4.9** Using Remark 4.4, the map $\varphi^2_{(4)} : H^2(\mathcal{M}^\bullet_{(4)}(C)) \to K^{(4)}_6(C)/K^{(2)}_2(F) \cup K^{(2)}_2(F) + K^{(3)}_5(k) \cup F^*_\mathbb{Q}$ in Corollary 4.8 can be lifted to a map $\varphi^2_{(4)} : H^2(\mathcal{M}^\bullet_{(4)}(C)) \to K^{(4)}_6(C)/K^{(2)}_4(F) \cup K^{(2)}_2(F)$ by changing $\varphi^2_{(4)}(\alpha)$ with elements in $K^{(3)}_5(k) \cup F^*_\mathbb{Q}$.

Large parts of the proof of Theorem 4.6 work for general $n$, and we give most of it in this context. However, although the method employed probably works for all $n$, the combinatorics at a certain stage get rather out of hand, so we restrict our attention to $n = 3$ at some point.

Starting with $\alpha = \sum_j c_j [f_j]_n \otimes g_j$ in $H^2(\mathcal{M}^\bullet_{(n+1)}(F))$ we begin with creating an element $\alpha_1$ in $K_{n+1}^{(n+1)}(X_{F,\text{loc}}^{n-1}; \square^{n-1})$. Let $\{A_1, \ldots, A_l\} \subset \{F^*\}$ be a basis of $<f_j, g_j> \subset F_\mathbb{Q}^*$ obtained by first choosing a basis of $<f_j>$ among the $f_j$'s and then extending to a basis of $<f_j, g_j>$. Write

$$f_j = \prod_k A_k^{s_{kj}} \qquad \text{and} \qquad g_j = \prod_k A_k^{t_{kj}} \tag{4.1}$$

in $F_\mathbb{Q}^*$. Let $F_j(t) = \dfrac{t - f_j}{t - 1} \prod \left( \dfrac{t - A_k}{t - 1} \right)^{-s_{kj}} \in (1 + I)^*$.

Let $J = (i_1 i_2 \cdots i_k)$ with all $i_j \in \{1, \ldots, n-1\}$ be a sequence of distinct elements, and let $J^{\text{ord}} = (j_1 j_2 \cdots j_k)$ be the ordered version of $J$, i.e., $J$ and $J^{\text{ord}}$ have the same elements, and $j_1 < j_2 < \ldots < j_k$. We shall write $(-1)^J$ for the sign of the permutation $\begin{pmatrix} i_1 & \cdots & i_k \\ j_1 & \cdots & j_k \end{pmatrix}$, and $(-1)^{j \in J}$ to mean $(-1)^l$ if $j = i_l$. If $J_1$ and $J_2$ are disjoint tuples, we write $J_1 J_2$ for there juxtaposition. If $j \in J$, write $J \setminus \{j\}$ for the $(|J| - 1)$–tuple obtained by deleting $j$ from $J$. Note that if $1, \ldots, n-1$ are all in $J$ then, by adding $j$ into the $j$–th position of $J \setminus \{j\}$, and moving it up front and then to its position in $J$, we find

$$(-1)^{J \setminus \{j\}} = (-1)^{j-1}(-1)^{(j)J \setminus \{j\}} = (-1)^j (-1)^{j \in J} (-1)^J. \tag{4.2}$$

For a set $I \subset \{1, \ldots, n-1\}$ we identify $I$ with the ordered tuple it defines by ordering its elements, and similarly for its complement $I_c = \{1, \ldots, n-1\} \setminus I$.

For $I = \{i_1 \ldots, i_k\}$ and $I_c = \{j_1, \ldots, j_{n-1-k}\}$ as above, let $F_j^I = F_j(t_{i_1}) \cup \cdots \cup F_j(t_{i_k})$ and $[f_j]_{|I_c|+1}^{I_c} = [f_j]_{n-k}$ seen as element of $K_{n-k}^{(n-k)}(X_\text{loc}^{n-1-k}; \square^{n-1-k})$ with coordinates $t_{j_1}, \ldots, t_{j_{n-1-k}}$. Then we let

$$d_i^A F_j^I \cup [f_j]_{n-1-k}^{I_c} \cup g_j$$
$$= \text{contribution to the boundary at } t_i = A_k\text{'s coming from the } F_j$$
$$= (-1)^{i \in II_c} \sum_k s_{kj} F_j^{I \setminus \{i\}} \cup [f_j]_{|I_c|+1}^{I_c} \cup g_{j|t_i = A_k}$$

(noting that $(-1^{i \in II_c})$ if $i \in I$) if $i \in I$ and zero otherwise. Similarly, we define

$$d_i^f F_j^I \cup [f_j]_{|I_c|+1}^{I_c} \cup g_j$$
$$= \text{contribution to the boundary at } t_i = f_j\text{'s coming from the } F_j$$
$$= - (-1)^{i \in II_c} F_j^{I \setminus \{i\}} \cup [f_j]_{|I_c|+1}^{I_c} \cup g_{j|t_i = f_j}$$

if $i \in I$ and zero otherwise, and

$$d_i^{[\,]} F_j^I \cup [f_j]_{|I_c|+1}^{I_c} \cup g_j$$
$$= \text{contribution to the boundary at } t_i = f_j\text{'s coming from the } [f]_{|I_c|}$$
$$= (-1)^{i \in II_c} F_j^I \cup [f_j]_{|I_c \setminus \{i\}|+1}^{I_c \setminus \{i\}} \cup g_{j|t_i = f_j},$$

(because $(-1)^{|I|}(-1)^{i \in I_c} = (-1)^{i \in II_c}$ if $i \in I_c$) if $i \notin I$ and zero otherwise. Note that $d_i = d_i^A + d_i^f + d_i^{[\ ]}$.

In the commutative diagram

$$\begin{array}{ccc} K_n^{(n+1)}(X^n; \square^n) & \xrightarrow{\cong} & K_{n+1}^{(n+1)}(X^{n-1}; \square^{n-1}) \\ \downarrow & & \downarrow \\ K_n^{(n+1)}(X \times X_{\text{loc}}^{n-1}; \square^n) & \xrightarrow{\cong} & K_{n+1}^{(n+1)}(X_{\text{loc}}^{n-1}; \square^{n-1}) \end{array}$$

we know by Theorem 2.4 that $\varphi_{(n+1)}^2$ maps $\sum_j c_j [f_j]_n \otimes g_j$ to $\pm \sum_j c_j [f_j]_n \cup g_j$ (modulo $(1+I)^* \tilde{\cup} K_n^{(n)}(X_{\text{loc}}^{n-2}; \square^{n-2})$) in $K_{n+1}^{(n+1)}(X_{\text{loc}}^{n-1}; \square^{n-1})$. The complex (from a spectral sequence analogous to (2.1))

$$K_{n+1}^{(n+1)}(X_{\text{loc}}^{n-1}; \square^{n-1}) \to \coprod K_n^{(n)}(X_{\text{loc}}^{n-2}; \square^{n-2}) \to \ldots \to K_2^{(2)}(F)$$

has the acyclic subcomplex

$$(1+I)^* \tilde{\cup} K_n^{(n)}(X_{\text{loc}}^{n-2}; \square^{n-2}) \to d(\ldots) + \coprod (1+I)^* \tilde{\cup} K_{n-1}^{(n-1)}(X_{\text{loc}}^{n-3}; \square^{n-3}) \to \cdots$$
$$\cdots \to d(\ldots) + \coprod (1+I)^* \tilde{\cup} K_2^{(2)}(F) \to d(\ldots)$$

with quotient complex

$$\frac{K_{n+1}^{(n+1)}(X_{\text{loc}}^{n-1}; \square^{n-1})}{(1+I)^* \tilde{\cup} K_n^{(n)}(X_{\text{loc}}^{n-2}; \square^{n-2})} \to \frac{K_n^{(n)}(X_{\text{loc}}^{n-2}; \square^{n-2})}{(1+I)^* \tilde{\cup} K_{n-1}^{(n-1)}(X_{\text{loc}}^{n-3}; \square^{n-3})} \otimes F_{\mathbb{Q}}^* \to \qquad (4.3)$$

(The proof that the subcomplex is acyclic is completely analogous to the proof of Lemma 3.7 in [dJ1] which is based on the fact that the subcomplex is closed under multiplication by $(1+I)^*$ and every element in the subcomplex contains at least one factor in $(1+I)^*$, see also [dJ1, Remark 3.10].) We know that $\alpha = \sum_j c_j [f_j]_n \otimes g_j$ is mapped to zero under the map in (4.3), and want to lift it back (uniquely because of the acyclicity of the subcomplex) to $\alpha_1$ in the kernel of $K_{n+1}^{(n+1)}(X_{\text{loc}}^{n-1}; \square^{n-1}) \to \coprod K_n^{(n)}(X_{\text{loc}}^{n-2}; \square^{n-2})$. Note that this lift is the restriction to $K_{n+1}^{(n+1)}(X_{\text{loc}}^{n-1}; \square^{n-1})$ of the image under $\varphi_{(n+1)}^2$ in $K_{n+1}^{(n+1)}(X_{\text{loc}}^{n-1}; \square^{n-1})$. The same proof works over some suitable Zariski open part of $C$.

**Proposition 4.10** If

$$\sum_j c_j [f_j]_n \cup g_j \in \frac{K_{n+1}^{(n+1)}(X_{\text{loc}}^{n-1}; \square^{n-1})}{(1+I)^* \tilde{\cup} K_n^{(n)}(X_{\text{loc}}^{n-2}; \square^{n-2})}$$

has trivial boundary $\sum_j c_j [f_j]_{n-1} \cup g_j \otimes f_j$ in

$$\frac{K_n^{(n)}(X_{\text{loc}}^{n-2}; \square^{n-2})}{(1+I)^* \tilde{\cup} K_{n-1}^{(n-1)}(X_{\text{loc}}^{n-3}; \square^{n-3})} \otimes F_{\mathbb{Q}}^*$$

(resp. $K_2^{(2)}(F) \otimes F_{\mathbb{Q}}^*$ for $n = 2$) then

$$\alpha_1 = \sum_{I \subset \{1,\ldots,n-1\}} (-1)^{II_c} \sum_j c_j F_j^I \cup [f_j]_{n-|I|}^{I_c} \cup g_j$$

(with the sum over all subsets $I$ of $\{1, \ldots, n-1\}$ seen as tuples in ascending order) in $K_{n+1}^{(n+1)}(X_{\text{loc}}^{n-1}; \square^{n-1})$ has zero boundary in $K_n^{(n)}(X_{\text{loc}}^{n-2}; \square^{n-2})$.

**Proof** The proof will be on induction on $n$. We need the following lemma.

**Lemma 4.11** With $f_j = \prod_k A_k^{s_{kj}}$ as before,

$$\sum_j c_j s_{kj} [f_j]_{n-1} \cup g_j = 0 \quad \text{in} \quad \frac{K_n^{(n)}(X_{\text{loc}}^{n-2}; \square^{n-2})}{(1+I)^* \tilde{\cup} K_{n-1}^{(n-1)}(X_{\text{loc}}^{n-3}; \square^{n-3})}$$

(resp. $K_2^{(2)}(F) \otimes F_{\mathbb{Q}}^*$ for $n = 2$).

**Proof** Write the boundary in $\dfrac{K_n^{(n)}(X_{\text{loc}}^{n-2}; \square^{n-2})}{(1+I)^* \tilde{\cup} K_{n-1}^{(n-1)}(X_{\text{loc}}^{n-3}; \square^{n-3})} \otimes F_{\mathbb{Q}}^*$ in terms of $\cdots \otimes A_k$'s, and collect terms, remembering that the $A_k$'s form a basis of $<f_j, g_j> \subset F_{\mathbb{Q}}^*$, so they are independent in $F_{\mathbb{Q}}^*$.

We now compute the boundary of $\alpha_1$, doing it for all $t_i$. For $n = 2$ one checks easily that $\alpha_1$ has boundary $-\sum_j c_j [f_{j,k}]_{c_j} s_{kj}(1-f_j)^{-1} \cup g_j = 0$ by writing it out in terms of our chosen bases for $<f_j, g_j> \subset F_{\mathbb{Q}}^*$. For the higher $n$'s, as $d_i = d_i^A + d_i^f + d_i^{[\,]}$, we get three contributions. We start with the $d_i^A$–component. We find

$$\sum_{\substack{I \subset \{1,\ldots,n-1\} \\ i \in I}} (-1)^{i \in II_c}(-1)^{II_c} \sum_{j,k} c_j s_{kj} F^{I \setminus \{i\}} \cup [f_j]_{|I_c|+1}^{I_c} \cup g_{j|t_i = A_k}$$

$$= (-1)^i \sum_{J \subset \{1,\ldots,\hat{i},\ldots,n-1\}} (-1)^{JJ_c} \sum_{j,k} c_j s_{kj} F_j^J \cup [f_j]_{|J_c|+1}^{J_c} \cup g_{j|t_i = A_k}$$

by letting $J = I \setminus \{i\}$, and taking $J_c$ in $\{1, \ldots, \hat{i}, \ldots, n-1\}$, and using (4.2). By induction on $n$ and Lemma 4.11, this equals zero.

For the contribution from $d_i^f$ we get

$$-\sum_{\substack{I \subset \{1,\ldots,n-1\} \\ i \in I}} (-1)^{i \in II_c}(-1)^{II_c} \sum_j c_j F^{I \setminus \{i\}} \cup [f_j]_{|I_c \setminus \{i\}|+1}^{I_c \setminus \{i\}} \cup g_{j|t_i = f_j}$$

$$= -(-1)^i \sum_{J \subset \{1,\ldots,\hat{i},\ldots,n-1\}} (-1)^{JJ_c} \sum_j c_j F_j^J \cup [f_j]_{|J_c|+1}^{J_c} \cup g_{j|t_i = A_k}$$

again by letting $J = I \setminus \{i\}$, and taking $J_c$ in $\{1, \ldots, \hat{i}, \ldots, n-1\}$, and using (4.2).

For $d_i^{[\,]}$ we get a contribution

$$\sum_{\substack{I \subset \{1,\ldots,n-1\} \\ i \notin I}} (-1)^{i \in II_c}(-1)^{II_c} \sum_j c_j F^I \cup [f_j]_{|I_c \setminus \{i\}|+1}^{I_c \setminus \{i\}} \cup g_{j|t_i = f_j}$$

$$= (-1)^i \sum_{J \subset \{1,\ldots,\hat{i},\ldots,n-1\}} (-1)^{JJ_c} \sum_j c_j F^J \cup [f_j]_{|J_c|+1}^{J_c} \cup g_{j|t_i = f_j}$$

by letting $J = I \subset \{1, \ldots, \hat{i}, \ldots, n-1\}$, and taking $J_c$ in $\{1, \ldots, \hat{i}, \ldots, n-1\}$ as before, and using (4.2) again. Obviously, the contributions of $d_i^f$ and $d_i^{[\,]}$ cancel.

Note that $k(x)$ is a number field and the regulator is injective (up to torsion) on its $K$–theory, we can compute the boundary at the level of Deligne cohomology, so we now turn towards the regulator level. Consider the following commutative diagram with vertical maps being the regulators (into deRham cohomology as it is equal to the Deligne cohomology in all cases considered).

$$\begin{array}{ccccc}
H^2(\widetilde{\mathcal{M}}^\bullet_{(n+1)}(F)) & & & & \\
\downarrow & & & & \\
K_n^{(n+1)}(X^n; \square^n) & \xrightarrow{\sim} & K_{n+1}^{(n+1)}(X^{n-1}; \square^{n-1}) & \longrightarrow & K_{n+1}^{(n+1)}(X_{F,\mathrm{loc}}^{n-1}; \square^{n-1}) \\
\downarrow & & \downarrow & & \downarrow \\
H_{\mathrm{dR}}^{n+1}(X^n; \square^n; \mathbb{R}(n+1))^+ & \xrightarrow{\sim} & H_{\mathrm{dR}}^n(X^{n-1}; \square^{n-1}; \mathbb{R}(n))^+ & \longrightarrow & H_{\mathrm{dR}}^n(X_{U,\mathrm{loc}}^{n-1}; \square^{n-1}; \mathbb{R}(n))^+
\end{array}$$

We shall determine the regulator of $\alpha_1$ in $H_{\mathrm{dR}}^n(X_{U,\mathrm{loc}}^{n-1}; \square^{n-1}; \mathbb{R}(n))^+$ and then lift it back to $H_{\mathrm{dR}}^n(X_U^{n-1}; \square^{n-1}; \mathbb{R}(n))^+$. We begin with computing the indeterminacy in the lift. Write $X_{\mathrm{loc}}$ for $X_{U,\mathrm{loc}} = X_U \setminus \{t = f_j\}$ with $U$ a suitable Zariski open part as before, $H^p(X_{\mathrm{loc}}^q; \mathbb{R}(s))$ for $H_{\mathrm{dR}}^p(X_{\mathrm{loc}}^q; \square^q; \mathbb{R}(s))$ and consider the spectral sequence

$$\begin{array}{llll}
H_{\mathrm{dR}}^n(X_{\mathrm{loc}}^{n-1}; \mathbb{R}(n)) & \coprod H_{\mathrm{dR}}^{n-1}(X_{\mathrm{loc}}^{n-2}; \mathbb{R}(n-1)) & \coprod H_{\mathrm{dR}}^{n-2}(X_{\mathrm{loc}}^{n-3}; \mathbb{R}(n-2)) & \\
H_{\mathrm{dR}}^{n-1}(X_{\mathrm{loc}}^{n-1}; \mathbb{R}(n)) & \coprod H_{\mathrm{dR}}^{n-2}(X_{\mathrm{loc}}^{n-2}; \mathbb{R}(n-1)) & \coprod H_{\mathrm{dR}}^{n-3}(X_{\mathrm{loc}}^{n-3}; \mathbb{R}(n-2)) & \quad (4.4) \\
H_{\mathrm{dR}}^{n-2}(X_{\mathrm{loc}}^{n-1}; \mathbb{R}(n)) & \coprod H_{\mathrm{dR}}^{n-3}(X_{\mathrm{loc}}^{n-2}; \mathbb{R}(n-1)) & \coprod H_{\mathrm{dR}}^{n-4}(X_{\mathrm{loc}}^{n-3}; \mathbb{R}(n-2)) &
\end{array}$$

converging to $H_{\mathrm{dR}}^*(X^{n-1}; \mathbb{R}(n))$.

**Lemma 4.12** $H_{\mathrm{dR}}^n(X_{U,\mathrm{loc}}^k; \square^k) = 0$ if $n < k$.

**Proof** For $k = 0$ this is obvious. For $k \geq 1$ we have a spectral sequence (with notation as in (4.4))

$$\begin{array}{lll}
H_{\mathrm{dR}}^n(X_{\mathrm{loc}}^k; \mathbb{R}(j)) & \coprod H_{\mathrm{dR}}^{n-1}(X_{\mathrm{loc}}^{k-1}; \mathbb{R}(j-1)) & \coprod H_{\mathrm{dR}}^{n-2}(X_{\mathrm{loc}}^{k-2}; \mathbb{R}(j-2)) \\
H_{\mathrm{dR}}^{n-1}(X_{\mathrm{loc}}^k; \mathbb{R}(j)) & \coprod H_{\mathrm{dR}}^{n-2}(X_{\mathrm{loc}}^{k-1}; \mathbb{R}(j-1)) & \coprod H_{\mathrm{dR}}^{n-3}(X_{\mathrm{loc}}^{k-2}; \mathbb{R}(j-2)) \\
H_{\mathrm{dR}}^{n-2}(X_{\mathrm{loc}}^k; \mathbb{R}(j)) & \coprod H_{\mathrm{dR}}^{n-3}(X_{\mathrm{loc}}^{k-1}; \mathbb{R}(j-1)) & \coprod H_{\mathrm{dR}}^{n-4}(X_{\mathrm{loc}}^{k-2}; \mathbb{R}(j-2))
\end{array}$$

converging to $H_{\mathrm{dR}}^*(X_U^k; \square^k; \mathbb{R}(j))$. The only terms contributing to $H_{\mathrm{dR}}^n(X_U^k; \square^k; \mathbb{R}(j))$ are $H_{\mathrm{dR}}^{n-2p}(X_{U,\mathrm{loc}}^{k-p}; \square^{k-p})$'s, which are zero by induction for $p \geq 1$. The boundaries leaving $H_{\mathrm{dR}}^n(X_{U,\mathrm{loc}}^k; \square^j)$ land in $H_{\mathrm{dR}}^{n-2p+1}(X_{U,\mathrm{loc}}^{k-p}; \square^{k-p}; \mathbb{R}(j-p))$'s for $p \geq 1$, which are also zero by induction. Therefore we get isomorphisms $H_{\mathrm{dR}}^n(X_{U,\mathrm{loc}}^k; \square^k; \mathbb{R}(j)) \cong H_{\mathrm{dR}}^n(X_U^k; \square^k; \mathbb{R}(j)) \cong H_{\mathrm{dR}}^{n-k}(U; \mathbb{R}(j)) = 0$.

By Lemma 4.12 there are only two terms contributing to $H_{\mathrm{dR}}^n(X^{n-1}; \square^{n-1}; \mathbb{R}(n-1))$ in (4.4) and we have a short exact sequence

$$0 \to E_2^+ \to H_{\mathrm{dR}}^n(X^{n-1}; \square^{n-1}; \mathbb{R}(n))^+ \to E_3^+ \to 0$$

with $E_3^+$ the +–part of the $E_\infty = E_3$ term at the $H^n_{dR}(X^{n-1}_{U,\text{loc}};\square^{n-1};\mathbb{R}(n))$–position and $E_2^+$ the +–part of the $E_\infty = E_2$ term at the $\coprod H^{n-2}_{dR}(X^{n-2}_{U,\text{loc}};\square^{n-2};\mathbb{R}(n-1))$ position. Because $H^n_{dR}(X^{n-1};\square^{n-1};\mathbb{R}(n))^+$ is alternating for the action of $S_{n-1}$ (acting on everything by permuting the coordinates), $E_2^+$ and $E_3^+$ are alternating as well. As we are looking at the regulator of an element coming from $X^{n-1}_U$ of course its survives in the spectral sequence to $E_\infty$ and we can consider its projection in $E_3^+$. We move on to determining the $E_2$–term. For this we introduce

$$R_c = <d\,i\arg\prod_j \left(\frac{t-f_j}{t-1}\right)^{n_j} \text{ such that } \prod_j f_j \in k^*>_\mathbb{R} \subset H^1_{dR}(X_{U,\text{loc}};\square;\mathbb{R}(1))^+$$

Note that an element in $R_c$ is determined completely by its residue at the $t - f_j$'s.

**Lemma 4.13** The map $H^0_{dR}(U;\mathbb{R}(0))_{|t=f} \to H^2_{dR}(X_U;\square;\mathbb{R}(1)) \cong H^1_{dR}(U;\mathbb{R}(1))$ maps 1 to $\pm d\,i\arg f$

**Proof** Consider the situation $U = \mathbb{G}_m$ and $f = S$ ($S$ the coordinate on $\mathbb{G}_m$). We then have the exact sequence in relative deRham cohomology

$$\ldots \to H^0(X_{\text{loc}}) \to H^0(U)^{\oplus 2} \to H^1(X_{\text{loc}};\square) \to H^1(X_{\text{loc}}) \to H^1(U)^{\oplus 2} \to \ldots$$

As $H^1(X_{\text{loc}}) = <d i\arg\frac{t-S}{t-1}>_\mathbb{R} \oplus <d i\arg S>_\mathbb{R}$, the last map in the above sequence is injective. From the corresponding sequence with $X$ instead of $X_{\text{loc}}$ one then gets that

$$H^1(X;\square) \cong H^1(X_{\text{loc}};\square)$$

as $H^0(X_{\text{loc}}) \cong H^0(X)$. Considering the localization sequence

$$H^1(X;\square)\tilde{\to}H^1(X_{\text{loc}};\square) \to H^0(U) \to H^2(X_U;\square) \to H^2(X_{\text{loc}};\square)$$

we see that $H^2(X_{U,\text{loc}};)\tilde{\to}H^2(X_{U,\text{loc}})$ because $H^2(\square) = 0$. From the commutative diagram

$$\begin{array}{ccc} H^2(X_U;\square) & \longrightarrow & H^2(X_U) \\ \downarrow & & \downarrow \\ H^2(X_{U,\text{loc}};\square) & \stackrel{\cong}{\longrightarrow} & H^2(U) \end{array}$$

we see that the map $H^2(X_U;\square) \to H^2(X_{U,\text{loc}};\square)$ is the zero map because $H^2(X_U) \cong H^2(U) = 0$. Hence $H^0(U)\tilde{\to}H^2(X_U;\square)\tilde{\to}H^1(U)$. All this works with cohomology with $\mathbb{Z}$–coefficients, which gives the statement for $\mathbb{G}_m$. By pulling back to our original $U$ via $f$ we get the corresponding statement for $f$ and $U$.

**Lemma 4.14**

$$H^n(X^n_{U,\text{loc}};\square^n;\mathbb{R}(n))^{\text{alt}} \cong \left(\bigoplus_{k=0}^n R_c^k \cup H^{n-k}(X^{n-k}_k;\square^{n-k};\mathbb{R}(n-k))\right)^{\text{alt}}$$

**Proof** For $n = 0$ this is obvious, or for $n = 1$ consider the localization sequence

$$0 \to H^1_{\mathrm{dR}}(X_U; \square; \mathbb{R}(1)) \to H^1_{\mathrm{dR}}(X_{U,\mathrm{loc}}; \square; \mathbb{R}(1)) \to$$
$$\coprod_j H^0_{\mathrm{dR}}(U; \mathbb{R}(0))_{|t=f_j} \xrightarrow{\varphi} H^2(X_U; \square; \mathbb{R}(1))$$

By Lemma 4.13, the composition $H^0_{\mathrm{dR}}(U; \mathbb{R}(0))_{|t=f_j} \to H^2_{\mathrm{dR}}(X_U; \square; \mathbb{R}(1)) \cong H^1_{\mathrm{dR}}(U; \mathbb{R}(1))$ maps 1 to $\pm \mathrm{d}\, i \arg f_j$. So $\coprod_j a_j$ is in the kernel of $\varphi$ if and only if $\sum_j a_j \mathrm{d}\, i \arg f_j = 0$, which means that $\coprod_j a_j$ is the image of $\sum_j a_j \mathrm{d} i \arg \frac{t - f_j}{t - 1}$ in $H^1_{\mathrm{dR}}(X_{U,\mathrm{loc}}; \square; \mathbb{R}(1))$. If we show this is in $R_c$, we are done. We have an exact sequence

$$0 \to k^* \to F^* \to \{\mathrm{d} i \arg f_j\}$$

as one sees by considering the residue versus the divisor map. Tensoring with $\mathbb{R}$ we get

$$0 \to k^* \otimes_\mathbb{Z} \mathbb{R} \to F^* \otimes_\mathbb{Z} \mathbb{R} \to <\mathrm{d} i \arg f_j>_\mathbb{R}$$

from which it follows that $\sum_j a_j \mathrm{d} i \arg \frac{t - f_j}{t - 1}$ is in $R_c$, e.g., by considering a $\mathbb{Q}$–basis of $\mathbb{R}$. Because $R_c$ injects into $\coprod_j H^0_{\mathrm{dR}}(U; \mathbb{R}(0))_{|t=f_j}$ under the residue, we get the statement for $n = 1$.

For $n \geq 2$, we use induction. We have a spectral sequence

$$\begin{array}{lll}
H^n_{\mathrm{dR}}(X^n_{\mathrm{loc}}; \mathbb{R}(n))^{\mathrm{alt}} & \left(\coprod H^{n-1}_{\mathrm{dR}}(X^{n-1}_{\mathrm{loc}}; \mathbb{R}(n-1))\right)^{\mathrm{alt}} & \left(\coprod H^{n-2}_{\mathrm{dR}}(X^{n-2}_{\mathrm{loc}}; \mathbb{R}(n-2))\right)^{\mathrm{alt}} \\
H^{n-1}_{\mathrm{dR}}(X^n_{\mathrm{loc}}; \mathbb{R}(n))^{\mathrm{alt}} & \left(\coprod H^{n-2}_{\mathrm{dR}}(X^{n-1}_{\mathrm{loc}}; \mathbb{R}(n-1))\right)^{\mathrm{alt}} & \left(\coprod H^{n-3}_{\mathrm{dR}}(X^{n-2}_{\mathrm{loc}}; \mathbb{R}(n-2))\right)^{\mathrm{alt}} \\
H^{n-2}_{\mathrm{dR}}(X^n_{\mathrm{loc}}; \mathbb{R}(n))^{\mathrm{alt}} & \left(\coprod H^{n-3}_{\mathrm{dR}}(X^{n-1}_{\mathrm{loc}}; \mathbb{R}(n-1))\right)^{\mathrm{alt}} & \left(\coprod H^{n-4}_{\mathrm{dR}}(X^{n-2}_{\mathrm{loc}}; \mathbb{R}(n-2))\right)^{\mathrm{alt}}
\end{array}$$

converging to $H^*_{\mathrm{dR}}(X^n_{\mathrm{loc}}; \square^n; \mathbb{R}(n))^{\mathrm{alt}} \cong H^*_{\mathrm{dR}}(X^n_{\mathrm{loc}}; \square^n; \mathbb{R}(n))$. Introducing the notation $H^k_{\mathrm{loc}}$ for $H^k_{\mathrm{dR}}(X^k_{U,\mathrm{loc}}; \square^k; \mathbb{R}(k))$, by Lemma 4.12 everything below the line

$$H^n_{\mathrm{loc}} \to \left(\coprod H^{n-1}_{\mathrm{loc}}\right)^{\mathrm{alt}} \to \left(\coprod H^{n-2}_{\mathrm{loc}}\right)^{\mathrm{alt}} \to \cdots \qquad (4.5)$$

vanishes. Write $H^k$ for $H^k_{\mathrm{dR}}(X^k_{U,\mathrm{loc}}; \square^k; \mathbb{R}(k))^{\mathrm{alt}}$. The subcomplex of the corresponding non–alternating row given by

$$\oplus_{k=1}^n R_c^k \tilde{\cup} H^{n-k} \to \mathrm{d}(\ldots) + \coprod \oplus_{k=1}^{n-1} R_c^k \tilde{\cup} H^{n-1-k} \to \mathrm{d}(\ldots) + \coprod \oplus_{k=1}^{n-2} R_c^k \tilde{\cup} H^{n-2-k} \to \quad (4.6)$$

is acyclic as in the proof of Lemma 3.7 of [dJ1] (because it is closed under multiplication by elements in $R_c$ and the boundary is injective on $R_c$), so hence is its alternating part. Taking the quotient of the complex in (4.5) and the alternating part of its subcomplex (4.6) yields by induction on $n$ the row

$$\left(\frac{H^n_{\mathrm{dR}}(X^n_{U,\mathrm{loc}}; \square^n; \mathbb{R}(n))}{\oplus_{k=1}^n R_c^k \tilde{\cup} H^{n-k}}\right)^{\mathrm{alt}} \to H^{n-1}_{\mathrm{dR}}(X^{n-1}_k; \square^{n-1}; \mathbb{R}(n-1)) \otimes_\mathbb{Q} F^*_\mathbb{Q}/k^*_\mathbb{Q} \to \cdots \qquad (4.7)$$

Obviously the last map is zero. Because the composition

$$H_{\mathrm{dR}}^{n-1}(X_k^{n-1}; \square^{n-1}; \mathbb{R}(n-1)) \otimes_{\mathbb{Q}} F_{\mathbb{Q}}^*/k_{\mathbb{Q}}^* \to H_{\mathrm{dR}}^{n+1}(X_U^n; \square^n; \mathbb{R}(n))^{\mathrm{alt}} \cong H^1(U; \mathbb{R}(n))$$

maps $c \wedge \mathrm{d}h(t_1) \wedge \ldots \wedge \mathrm{d}h(t_{n-1}) \otimes f_j$ to $c \mathrm{d} i \arg f_j$ this is an injection. So the first map in (4.7) must also be zero, giving an identification

$$H_{\mathrm{dR}}^n(X_k^n; \square^n; \mathbb{R}(k))^{\mathrm{alt}} \cong H_{\mathrm{dR}}^n(X_U^n; \square^n; \mathbb{R}(k))^{\mathrm{alt}} \cong \left( \frac{H_{\mathrm{dR}}^n(X_{U,\mathrm{loc}}^n; \square^n; \mathbb{R}(n))}{\oplus_{k=1}^n R_c^k \cup H^{n-k}} \right)^{\mathrm{alt}}$$

from which the result is immediate because $\oplus_{k=1}^n R_c^k \cup H^{n-k}$ injects under the residue into $\coprod H^{n-1}$.

**Remark 4.15** The proof shows that the $H_{\mathrm{dR}}^{n-2}(X_{\mathbb{C}}^{n-1}; \square^{n-1}; \mathbb{R}(n-1))_{|t_i=f_j}$ generate $E_2 = E_2^{\mathrm{alt}}$. Hence by Lemma 4.13, a lift from $E_3 = E_3^{\mathrm{alt}}$ to $H_{\mathrm{dR}}^n(X_U^{n-1}; \square^{n-1}; \mathbb{R}(n))^{\mathrm{alt}} \cong H_{\mathrm{dR}}^1 U; \mathbb{R}(n)$ is determined up to $\mathbb{R}(n-1) <\mathrm{d} i \arg f_j>_{\mathbb{R}} \wedge \mathrm{d}h(t_1) \wedge \ldots \wedge \mathrm{d}h(t_{n-1})$ corresponding to $\mathbb{R}(n-1) <\mathrm{d} i \arg f_j>_{\mathbb{R}}$ under this isomorphism.

Note that we can compute the residue as follows. There is a commutative diagram

$$\begin{array}{ccc}
K_{n+1}^{(n+1)}(X_U^{n-1}; \square^{n-1}) & \longrightarrow & K_n^{(n)}(X_{k(x)}^{n-1}; \square^{n-1}) \\
\downarrow{\scriptstyle \mathrm{reg}} & & \downarrow{\scriptstyle \mathrm{reg}} \\
H_{\mathrm{dR}}^n(X_U^{n-1}; \square^{n-1}; \mathbb{R}(n)) & \longrightarrow & H_{\mathrm{dR}}^{n-1}(X_{k(x)}^{n-1}; \square^{n-1}; \mathbb{R}(n)) \\
\downarrow{\scriptstyle \sim} & & \downarrow{\scriptstyle \sim} \\
H_{\mathrm{dR}}^1(U; \mathbb{R}(n)) & \longrightarrow & H_{\mathrm{dR}}^0(k(x); \mathbb{R}(n-1))
\end{array}$$

If $\psi$ is in $H_{\mathrm{dR}}^1(U; \mathbb{R}(n))$ and $x$ is a point not in $U$, then $\mathrm{res}_x(\psi)$ in $H_{\mathrm{dR}}^0(k(x); \mathbb{R}(n-1))$ is given by

$$\pm \frac{1}{(2\pi i)^n} \int_{X^{n-1} \times S_x^1} \psi \wedge \mathrm{d}h(t_1) \wedge \ldots \wedge \mathrm{d}h(t_{n-1}) \wedge \mathrm{d} i \arg t_1 \wedge \ldots \wedge \mathrm{d} i \arg t_{n-1}$$

for $S_x^1$ a circle around $x$.

We can also replace $U$ with the closed set of $C$ by leaving out small (open) discs around the point $x \notin U$, without changing either the cohomology groups involved or the values of the integrals. We shall assume that from now on, so in particular $U$ is compact.

**Lemma 4.16** Suppose $\psi_1$ and $\psi_2$ in $H_{\mathrm{dR}}^{n+1}(X_U^{n-1}; \square^{n-1}; \mathbb{R}(n-1))$ satisfy condition (2.5). Then with $\overline{\omega} = \mathrm{d} \log \overline{t_1} \wedge \ldots \wedge \mathrm{d} \log \overline{t_{n-1}}$, we have an equality

$$\int_{X^{n-1} \times S_x^1} \psi_1 \wedge \overline{\omega} = \int_{X^{n-1} \times S_x^1} \psi_2 \wedge \overline{\omega}$$

The same holds if we replace $\overline{\omega}$ with $\mathrm{d} i \arg t_1 \wedge \ldots \wedge \mathrm{d} i \arg t_{n-1}$.

**Proof** It was proved in the proof of Proposition 4.6 of [dJ2] that $\psi_1 - \psi_2 = d\gamma$ where $\gamma$ satifies the conditions in (2.5) on a suitable blowup of $(\mathbb{P}_C^1)^{n-1}$, isomorphic to this over a suitable Zariski open part of $C$. With that, one checks easily using integration in each fibre, that $\int_{X^{n-1} \times S_x^1} d\gamma \wedge \overline{\omega} = 0$ as the holomorphic form has a zero along $t_i = 1$ for every $i$. Hence the result follows from Stokes' theorem.

**Remark 4.17** Note that if we represent the image of $E_2^+$ inside $H_{dR}^n(X_U^{n-1}; \square^{n-1}; \mathbb{R}(n))^+$ by forms given by $H_{dR}^0(U; \mathbb{R}(n-1))^+ \cdot <d\,i \arg f_{j\mathbb{R}}> \wedge dh(t_1) \wedge \ldots \wedge dh(t_{n-1})$ then

$$\coprod_{x \notin U} \int_{X^{n-1} \times S_x^1} \psi \wedge d\,i \arg t_1 \wedge \ldots \wedge d\,i \arg t_{n-1}$$

converges and maps $E_2^+$ to $H_{dR}^0(U; \mathbb{R}(n-1))^+ \cdot \text{res}<d\,i \arg f_j>$ in $\coprod_{x \notin U} \mathbb{R}(n-1)$.

**Proof** The computation of $E_2^+$ was carried out in the proof of Lemma 4.14. It is generated by the pushforward of the alternating version of $H_{dR}^{n-2}(X_k^{n-1}; \square^{n-1}; \mathbb{R}(n-1))_{|t_i = f_j}$. The rest is just a matter of integration and Lemma 4.13.

**Remark 4.18** Note that by Borel's theorem the image of $H_{dR}^0(U; \mathbb{R}(n-1))^+ \cdot \text{res}<d\,i \arg f_j>$ is exactly $\coprod_{x \notin U} \int_{S^1} \text{reg}(K_{2n-1}^{(n)}(k))<d\,i \arg f_j>_\mathbb{R}$

We now have the regulator $\text{reg}(\alpha_1)$ of $\alpha_1$ in $H_{dR}^n(X_{U,\text{loc}}^{n-1}; \square^{n-1}; \mathbb{R}(n))$. If we lift it back to $\psi \in H_{dR}^n(X_U^{n-1}; \square^{n-1}; \mathbb{R}(n))$ satifying (2.5) then by Remark 4.17 and Lemma 4.16

$$-\frac{1}{(2\pi i)^3} \coprod_{x \notin U} \int_{X^{n-1} \times S_x^1} \psi \wedge d\,i \arg t_1 \wedge \ldots \wedge d\,i \arg t_{n-1} \tag{4.8}$$

differs from the boundary

$$-\frac{1}{(2\pi i)^3} \coprod_{x \notin U} \int_{X^{n-1} \times S_x^1} \text{reg}(\alpha) \wedge d\,i \arg t_1 \wedge \ldots \wedge d\,i \arg t_{n-1} = \tag{4.9}$$

by an element in

$$\frac{1}{2\pi i} \coprod_{x \notin U} \int_{S_x^1} \text{reg}(K_{2n-1}^{(n)}(k))<d\,i \arg f_j>_\mathbb{R}.$$

Because we shall see that all values in (4.8) and necessarily (4.9) are in the image of $K_{2n-1}^{(n)}(k(x))$, it follows that $\psi = \text{reg}(\alpha) + \text{reg}(\beta)$ for some $\beta$ in $<f_j> \cup K_{2n-1}^{(n)}(k)$. From this we get Theorem 4.6.

We now turn our attention to the explicit lift of $\text{reg}(\alpha_1)$, starting out in general, but specializing to the case $n = 3$ at some stage.

**Lemma 4.19** Suppose we have $(n-1-k)$–forms $\varphi_k(t_{k+1}, \ldots, t_{n-1})$ for $k = 0, \ldots, n-1$, such that $\varphi_{k-1}(0, t_{k+1}, \ldots, t_{n-1}) = d\varphi_k(t_{k+1}, \ldots, t_{n-1})$ for $k = 1, \ldots, n-1$. Assume moreover that each $\varphi_k$ is alternating for the action of $S_{n-k-1}$, and $\varphi_k$ vanishes for $t_j = \infty$, $j = k+1, \ldots, n-1$. Let $\rho(t)$ be a bump form around $t = 0$. Then the form

$$\sum_{I \subset \{1, \ldots, n-1\}} (-1)^{|I|} (-1)^{II_c} D^I \varphi_{|I|}(t^{I_c})$$

where the sum is over all ordered tuples $I$ of $\{1, \ldots, n-1\}$, (and its complement $I_c$ is given the ascending ordering), $t^{I_c} = t_{j_1}, \ldots, t_{j_k}$ and if $I_c = \{j_1, \ldots, j_k\}$ with $j_1 < \ldots < j_k$, and $D^J \varphi = \mathrm{d}\left(\rho(t_{j_1})\mathrm{d}\left(\rho(t_{j_2}) \cdots \mathrm{d}\left(\rho(t_{j_k})\varphi\right)\cdots\right)\right)$ if $J = (j_1 \cdots j_k)$, is an alternating form vanishing at $t_j = 0$ for $k = 1, \ldots, n-1$.

**Proof** For $\sigma$ in $S_{n-1}$, we have that

$$\begin{aligned}
\sigma^*\left((-1)^{II_c} D^I \varphi_{|I|}(t^{I_c})\right) &= (-1)^{II_c} D^{\sigma(I)} \varphi_{|\sigma(I)|}(t^{\sigma(I_c)}) \\
&= (-1)^\sigma (-1)^{\sigma(I)\sigma(I_c)} D^{\sigma(I)} \varphi_{|\sigma(I)|}(t^{\sigma(I_c)}) \\
&= (-1)^\sigma (-1)^{\sigma(I)(\sigma(I))_c} D^{\sigma(I)} \varphi_{|\sigma(I)|}(t^{\sigma(I)_c})
\end{aligned}$$

by replacing $\sigma(I_c)$ with $(\sigma(I))_c$, i.e., ordering it. This shows the form is alternating. Note that $\left(D^I \varphi\right)_{|t_j=0} = 0$ unless $j \notin I$ or is its last element. Therefore, when restricting to $t_j = 0$ we get

$$\sum_{\substack{J \subset \{1, \ldots, \hat{j}, \ldots, n-1\} \\ I = J(j)}} (-1)^{|J|+1} (-1)^{J(j) J_c^{\mathrm{ord}}} D^J \mathrm{d}\varphi_{|J|+1}(t^{J_c^{\mathrm{ord}}})$$

$$+ \sum_{\substack{J \subset \{1, \ldots, \hat{j}, \ldots, n-1\} \\ I = J}} (-1)^{|J|} (-1)^{J((j)J_c)^{\mathrm{ord}}} D^J \varphi_{|J|}(t^{((j)J_c)^{\mathrm{ord}}})$$

with $t_j = 0$, and $J_c^{\mathrm{ord}}$ taken in $\{1, \ldots, \hat{j}, \ldots, n-1\}$. Using (4.2) and the conditions on the $\varphi_k$'s we have

$$\begin{aligned}
(-1)^{J(j)J_c^{\mathrm{ord}}} &= (-1)^j (-1)^{|J|+1} (-1)^{JJ_c^{\mathrm{ord}}} \\
\mathrm{d}\varphi_{|J|+1}(t^{J_c^{\mathrm{ord}}}) &= \varphi_{|J|}(0, t^{J_c^{\mathrm{ord}}}) \\
\varphi_{|J|}(t^{((j)J_c)^{\mathrm{ord}}})_{|t_j=0} &= -(-1)^{j \in ((j)J_c)^{\mathrm{ord}}} \varphi_{|J|}(0, t^{J_c^{\mathrm{ord}}}) \\
(-1)^{J((j)J_c)^{\mathrm{ord}}} &= (-1)^j (-1)^{j \in ((j)J_c)^{\mathrm{ord}}} (-1)^{|J|} (-1)^{JJ_c^{\mathrm{ord}}}
\end{aligned}$$

so everything cancels.

It turns out that for applying Lemma 4.19 with $\varphi_0 = \varepsilon_n$, writing down the forms involved is quite messy. We therefore assume from now on that $n = 3$. (The case $n = 2$ was done before in [dJ2].) In this case we shall carry out the lift explicitly.

We need some identities between forms (all $f$'s, $g$'s are functions on $C$).

$$\begin{aligned}
\mathrm{d}i \arg f_1 \wedge \mathrm{d}i \arg f_2 &= -\mathrm{d}\log|f_1| \wedge \mathrm{d}\log|f_2| \\
\mathrm{d}i \arg f_1 \wedge \mathrm{d}\log|f_2| &= -\mathrm{d}\log|f_1| \wedge \mathrm{d}i \arg f_2
\end{aligned}$$

Both identities follow easily by considering the $\pi_0$ or $\pi_1$ of $\mathrm{d}\log f_1 \wedge \mathrm{d}\log f_2 = 0$. We let

$$\sigma(f_1, f_2) = \log|f_1|\mathrm{d}i \arg f_2 - \log|f_2|\mathrm{d}i \arg f_1,$$

so $\mathrm{d}\sigma(f_1, f_2) = 0$.

Let $U = \mathbb{G}_m \setminus \{1\}$, and let $X_{U,\mathrm{loc}} = X_U \setminus \{t = S\}$. We want to write down the explicit elements $(\omega_n, \varepsilon_n) \in H_\mathcal{D}^n(X_{\mathrm{loc},U}^{n-1}; \square^{n-1}; \mathbb{R}(n))^+$ that are the images of $[S]_n$ under the regulator for $n = 1, 2$ and $3$, see Section 2.

For $n = 1$, we have $\tilde{\varepsilon}_1 = \varepsilon_1 = \log|1 - f|$ and $\omega_1 = \mathrm{d}\log(1-f)$. For $n = 2$, we have

$$\omega_2 = -\mathrm{d}\log\frac{t-f}{t-1} \wedge \mathrm{d}\log(1-f)$$

To find $\varepsilon_2$, let

$$\tilde{\varepsilon}_2 = \log|1-f|\mathrm{d}i\arg\left(\frac{t-f}{t-1}\right) - \log\left|\frac{t-f}{t-1}\right|\mathrm{d}i\arg(1-f)$$

and if we specialize this to $t = 0$, we find this equals

$$\mathrm{d}\varphi_1^{(2)}$$

with $\varphi_1^{(2)} = -P_{2,\mathrm{Zag}}(f)$. Then $\varepsilon_2 = \tilde{\varepsilon}_2 - \mathrm{d}\left(\rho(t)\varphi_1^{(2)}\right)$ with $\rho(t)$ a bump form around $t = 0$. (This is the correct $\varepsilon_2$, see [dJ2].)

Finally, for $n = 3$, we have that

$$\omega_3 = \mathrm{d}\log\frac{t_1-f}{t_1-1} \wedge \mathrm{d}\log\frac{t_2-f}{t_2-1} \wedge \mathrm{d}\log(1-f).$$

In order to find $\varepsilon_3$, let

$$\tilde{\varepsilon}_3 = \log|1-f|\mathrm{d}i\arg\left(\frac{t_1-f}{t_1-1}\right) \wedge \mathrm{d}i\arg\left(\frac{t_2-f}{t_2-1}\right)$$
$$+ \frac{2}{3!}\log|1-f|\mathrm{d}\log\left|\frac{t_1-f}{t_1-1}\right| \wedge \mathrm{d}\log\left|\frac{t_2-f}{t_2-1}\right|$$
$$- \log\left|\frac{t_1-f}{t_1-1}\right|\mathrm{d}i\arg(1-f) \wedge \mathrm{d}i\arg\left(\frac{t_2-f}{t_2-1}\right)$$
$$- \frac{2}{3!}\log\left|\frac{t_1-f}{t_1-1}\right|\mathrm{d}\log|1-f| \wedge \mathrm{d}\log\left|\frac{t_2-f}{t_2-1}\right|$$
$$+ \log\left|\frac{t_2-f}{t_2-1}\right|\mathrm{d}i\arg(1-f) \wedge \mathrm{d}i\arg\left(\frac{t_1-f}{t_1-1}\right)$$
$$+ \frac{2}{3!}\log\left|\frac{t_2-f}{t_2-1}\right|\mathrm{d}\log|1-f| \wedge \mathrm{d}\log\left|\frac{t_1-f}{t_1-1}\right|.$$

Specializing to $t_1 = 0$ we find after some computation that we get $\mathrm{d}\varphi_1^{(3)}$ where

$$\varphi_1^{(3)}(t_2) = -P_{2,\mathrm{Zag}}(f)\mathrm{d}i\arg\left(\frac{t_2-f}{t_2-1}\right) - \frac{2}{3}\log|1-f|\log\left|\frac{t_2-f}{t_2-1}\right|\mathrm{d}\log|f|$$
$$- \frac{1}{3}\log|1-f|\log|f|\mathrm{d}\log\left|\frac{t_2-f}{t_2-1}\right|.$$

Finally, putting $t_2 = 0$ in $\varphi_1^{(3)}$ we find that we get $\mathrm{d}\varphi_2^{(3)}$ with

$$\varphi_2(3)(f) = -P_{3,\mathrm{Zag}}(f) - \frac{1}{2}\log^2|f|\log|1-f|.$$

Putting everything together as in Lemma 4.19, we put

$$\varepsilon_3 = \tilde{\varepsilon}_3 - d[\rho(t_1)\varphi_1^{(3)}(t_2)] + d[\rho(t_2)\varphi_1^{(3)}(t_1)] + d[\rho(t_1)d[\rho(t_2)\varphi_2^{(3)}(f)]] - d[\rho(t_2)d[\rho(t_1)\varphi_2^{(3)}(f)]].$$

We check that $(\omega_3, \varepsilon_3)$ is the class of the regulator of $[f]_3$. Let $(\omega_3, \varepsilon_3')$ be the regulator. Then $\varepsilon_3 - \varepsilon_3' \in H^2_{dR}(X^2_{U,\text{loc}}; \square^2; \mathbb{R}(2))^+$. By (4.4) and Lemma 4.12, $H^2_{dR}(X^2; \square^2; \mathbb{R}(2))^+$ can be computed as the kernel of the mpa

$$H^2_{dR}(X^2_{U,\text{loc}}; \square^2; \mathbb{R}(2))^{+,\text{alt}} \xrightarrow{\text{res}} \left(\coprod H^1_{dR}(X^1_{U,\text{loc}}; \square^1; \mathbb{R}(1))^+\right)^{\text{alt}}$$

so if we can show that both $\varepsilon_3$ and $\varepsilon_3'$ have the same residue $\varepsilon_{2|t_1=S} - \varepsilon_{2|t_2=S}$ then they differ by an element of $H^2_{dR}(X^2; \square^2; \mathbb{R}(2))^+ \cong \mathbb{R}(2)$, and we can check that they are the same by specializing to a fixed value and integrating. From the exact sequence

$$0 \to H^1_{dR}(X_U; \square; \mathbb{R}(1))^+ \to H^1_{dR}(X^1_{U,\text{loc}}; \square^1; \mathbb{R}(1))^+ \to H^0_{dR}(U; \mathbb{R}(0))^+$$

we see, as $H^1_{dR}(X_U; \square; \mathbb{R}(1))^+ \cong H^0_{dR}(U; \mathbb{R}(1))^+ = 0$ that we do not lose any information about the residue in $H^1_{dR}(X^1_{U,\text{loc}}; \square^1; \mathbb{R}(1))^+$ by specializing $S$ to a constant. Therefore, assuming $S = c$ is constant in $\mathbb{Q}$ such that $\rho(c) = 0$ from now on, we find at $t_1 = c$ that the residue is

$$\log|1-c|di\arg\left(\frac{t_2-c}{t_2-1}\right)2 + d(P_{2,\text{Zag}}(c)\rho_1(t_2)) = \varepsilon_2(t_2,c)$$

as desired. So $\varepsilon_3 - \varepsilon_3'$ is an element of $H^2_{dR}(X^2; \square^2; \mathbb{R}(2))^+ \cong \mathbb{R}(2)$. Note that again in order to check that they are identical, we can specialize to $S = c$, so that we can check that they are the same class by computing

$$\int_{X^2} \varepsilon_3(t_1, t_2, c) \wedge di\arg t_1 \wedge di\arg t_2 = (2\pi i)^2 2 P_{3,\text{Zag}}(c).$$

Because this is the answer for $\varepsilon_3$ (see [dJ1, Proposition 4.1] with the correct sign), we conclude that $\varepsilon_3' = \varepsilon_3$.

Recall that we had fixed an element $\alpha = \sum_j c_j[f_j]_3 \otimes g_j$ in $H^2(\mathcal{M}^\bullet_{(4)}(F))$. Before writing down the corresponding regulator, we deduce some identities. We have the identity

$$\sum_j c_j[f_j]_2 \otimes (f_j \wedge g_j) = 0. \tag{4.10}$$

By applying $d \otimes id$ to it and writing it with respect to our basis of $<f_j, g_j>$ we find that for all $k$

$$\sum_j c_j s_{kj}(1-f_j) \otimes (f_j \otimes g_j - g_j \otimes f_j) = 0. \tag{4.11}$$

By using our basis once more we find that for all $k$ and $l$

$$\sum_j c_j s_{kj} s_{lj}(1-f_j) \otimes g_j = \sum_j c_j s_{kj} t_{lj}(1-f_j) \otimes f_j \tag{4.12}$$

We shall also need

$$\sum_{j,k} c_j s_{kj}[f_j]_2 \otimes g_j = \sum_{j,k,l} c_j s_{kj} t_{lj}[f_j]_2 \otimes A_l = \sum_{j,k,l} c_j s_{lj} t_{kj}[f_j]_2 \otimes A_l = \sum_{j,k} c_j[f_j]_2 \otimes f_j = 0$$

in $H^2(\mathcal{M}_{(3)}^\bullet(F))$, where we used (4.10) and used our basis again. Similarly we see from (4.10) that for all $k$ and $l$

$$\sum_j c_j s_{kj} t_{lj} [f_j]_2 = \sum_j c_j s_{lj} t_{kj} [f_j]_2. \tag{4.13}$$

Note also that the map $f \mapsto P_{2,\text{Zag}}(f)$ factors though $\widetilde{M}_{(2)}(F)$, so that an element of $\widetilde{M}_{(2)}(F)$ gives rise to a continuous function on $C$, differentiable where $f$ has no zeros or poles or assumes the value 1.

From Proposition 4.10 we obtain the element

$$\sum_j c_j [f_j]_3 \cup g_j$$
$$+ \sum_j c_j F_j(t_1) \cup [f_j]_2(t_2) \cup g_j$$
$$- \sum_j c_j F_j(t_2) \cup [f_j]_2(t_1) \cup g_j$$
$$+ \sum_j c_j F_j(t_1) \cup F_j(t_2) \cup [f_j]_1 \cup g_j$$

in $K_4^{(4)}(X_{\text{loc}}^2; \square^2)$. (Recall that $F_j(t) = \dfrac{t - f_j}{t - 1} \prod \left( \dfrac{t - A_k}{t - 1} \right)^{-s_{kj}} \in (1 + I)^*$ if $f_j = \prod_k A_k^{s_{kj}}$ as in (4.1).) According to (2.3), it has regulator

$$\begin{aligned}
& \sum_j c_j \varepsilon_3(t_1, t_2, f_j) \wedge \mathrm{d}i \arg g_j - \log|g_j| \pi_3 \omega_3(t_1, t_2, f_j) \\
& - \sum_j c_j \mathrm{d}i \arg F_j(t_1) \wedge (\varepsilon_2(t_2, f_j) \wedge \mathrm{d}i \arg g_j + \log|g_j| \pi_2 \omega_2(t_2, f_j)) \\
& + \sum_j c_j \mathrm{d}i \arg F_j(t_2) \wedge (\varepsilon_2(t_1, f_j) \wedge \mathrm{d}i \arg g_j + \log|g_j| \pi_2 \omega_2(t_1, f_j)) \\
& + \sum_j c_j \mathrm{d}i \arg F_j(t_1) \wedge \mathrm{d}i \arg F_j(t_2) \wedge \sigma(1 - f_j, g_j)
\end{aligned} \tag{4.14}$$

in $H^3_{\text{dR}}(X_{U,\text{loc}}^3; \square^3; \mathbb{R}(3))^+$. We want to lift it back to $H^3_{\text{dR}}(X_U^3; \square^3; \mathbb{R}(3))^+$. After some rewriting, one finds that $\tilde{\varepsilon}_3 \wedge \mathrm{d}i \arg g - \log|g| \pi_3 \omega_3$ equals

$$\mathrm{d}i \arg \left( \frac{t_1 - f}{t_1 - 1} \right) \wedge \mathrm{d}i \arg \left( \frac{t_2 - f}{t_2 - 1} \right) \wedge \sigma(1 - f, g)$$
$$+ \mathrm{d}\left[ \left( \frac{1}{3} \log|g| \mathrm{d}i \arg(1 - f) + \frac{1}{6} \sigma(g, 1 - f) \right) \wedge \right.$$
$$\left( \log\left| \frac{t_1 - f_j}{t_1 - 1} \right| \mathrm{d}\log\left| \frac{t_2 - f_j}{t_2 - 1} \right| - \log\left| \frac{t_2 - f_j}{t_2 - 1} \right| \mathrm{d}\log\left| \frac{t_1 - f_j}{t_1 - 1} \right| \right)$$
$$+ \log|g| \log\left| \frac{t_1 - f}{t_1 - 1} \right| \mathrm{d}\log|1 - f| \wedge \mathrm{d}i \arg\left( \frac{t_2 - f}{t_2 - 1} \right)$$
$$\left. - \log|g| \log\left| \frac{t_2 - f}{t_2 - 1} \right| \mathrm{d}\log|1 - f| \wedge \mathrm{d}i \arg\left( \frac{t_1 - f}{t_1 - 1} \right) \right].$$

For

$$\mathrm{d}i \arg F(t_1) \wedge (\tilde{\varepsilon}_2(t_2, f) \wedge \mathrm{d}i \arg g + \log|g| \pi_2 \omega_2(t_2, f))$$

we obtain
$$d i \arg F(t_1) \wedge d i \arg \left(\frac{t_2 - f}{t_2 - 1}\right) \wedge \sigma(1 - f, g)$$
$$+ d\left[-\log|g| \log\left|\frac{t_2 - f}{t_2 - 1}\right| d \log|1 - f| \wedge d i \arg F(t_1)\right].$$

Putting all this together, we get that the form in (4.14) after replacing all $\varepsilon$'s with $\tilde{\varepsilon}$'s, is given by

$$\begin{aligned}
\text{reg}^{\sim} = \sum_j c_j d \Bigg(&\frac{1}{6} \log\left|\frac{t_1 - f_j}{t_1 - 1}\right| \sigma(g_j, 1 - f_j) \wedge d \log\left|\frac{t_2 - f_j}{t_2 - 1}\right| \\
&- \frac{1}{6} \log\left|\frac{t_2 - f_j}{t_2 - 1}\right| \sigma(g_j, 1 - f_j) \wedge d \log\left|\frac{t_1 - f_j}{t_1 - 1}\right| \\
&+ \log|g_j| \log\left|\frac{t_1 - f_j}{t_1 - 1}\right| d \log|1 - f_j| \wedge d i \arg\left(\frac{t_2 - f_j}{t_2 - 1}\right) \\
&+ \frac{1}{3} \log|g_j| \log\left|\frac{t_1 - f_j}{t_1 - 1}\right| d i \arg(1 - f_j) \wedge d \log\left|\frac{t_2 - f_j}{t_2 - 1}\right| \\
&- \log|g_j| \log\left|\frac{t_2 - f_j}{t_2 - 1}\right| d \log|1 - f_j| \wedge d i \arg\left(\frac{t_1 - f_j}{t_1 - 1}\right) \\
&- \frac{1}{3} \log|g_j| \log\left|\frac{t_2 - f_j}{t_2 - 1}\right| d i \arg(1 - f_j) \wedge d \log\left|\frac{t_1 - f_j}{t_1 - 1}\right| \\
&+ \log|g_j| \log\left|\frac{t_2 - f_j}{t_2 - 1}\right| d \log|1 - f_j| \wedge d i \arg F_j(t_1) \\
&- \log|g_j| \log\left|\frac{t_1 - f_j}{t_1 - 1}\right| d \log|1 - f_j| \wedge d i \arg F_j(t_2) \Bigg) \\
&+ \sum_{j,k,l} c_j s_{kj} s_{lj} d i \arg\left(\frac{t_1 - A_k}{t_1 - 1}\right) \wedge d i \arg\left(\frac{t_2 - A_l}{t_2 - 1}\right) \wedge \sigma(1 - f_j, g_j)
\end{aligned}$$

Using (4.12) the last part can be written as

$$d \left(\sum_{j,k,l} c_j s_{kj} t_{lj} d i \arg\left(\frac{t_1 - A_k}{t_1 - 1}\right) \wedge d i \arg\left(\frac{t_2 - A_l}{t_2 - 1}\right) \wedge \sigma(1 - f_j, f_j)\right)$$
$$= \sum_{j,k,l} c_j s_{kj} t_{lj} d \left(-P_{2,\text{Zag}}(f_j) d i \arg\left(\frac{t_1 - A_k}{t_1 - 1}\right) \wedge d i \arg\left(\frac{t_2 - A_l}{t_2 - 1}\right)\right)$$

because $\sigma(1 - f, f) = -d P_{2,\text{Zag}}(f)$. Therefore $\text{reg}^{\sim}$ above equals $d \psi_1(t_1, t_2)$ with $\psi_1(t_1, t_2)$ given by

$$\sum_j c_j \frac{1}{6} \log\left|\frac{t_1 - f_j}{t_1 - 1}\right| \sigma(g_j, 1 - f_j) \wedge d \log\left|\frac{t_2 - f_j}{t_2 - 1}\right|$$
$$- \frac{1}{6} \log\left|\frac{t_2 - f_j}{t_2 - 1}\right| \sigma(g_j, 1 - f_j) \wedge d \log\left|\frac{t_1 - f_j}{t_1 - 1}\right|$$

$$+\log|g_j|\log\left|\frac{t_1-f_j}{t_1-1}\right|\operatorname{d}\log|1-f_j|\wedge\operatorname{d}i\arg\left(\frac{t_2-f_j}{t_2-1}\right)$$

$$+\frac{1}{3}\log|g_j|\log\left|\frac{t_1-f_j}{t_1-1}\right|\operatorname{d}i\arg(1-f_j)\wedge\operatorname{d}\log\left|\frac{t_2-f_j}{t_2-1}\right|$$

$$-\log|g_j|\log\left|\frac{t_2-f_j}{t_2-1}\right|\operatorname{d}\log|1-f_j|\wedge\operatorname{d}i\arg\left(\frac{t_1-f_j}{t_1-1}\right)$$

$$-\frac{1}{3}\log|g_j|\log\left|\frac{t_2-f_j}{t_2-1}\right|\operatorname{d}i\arg(1-f_j)\wedge\operatorname{d}\log\left|\frac{t_1-f_j}{t_1-1}\right|$$

$$+\log|g_j|\log\left|\frac{t_2-f_j}{t_2-1}\right|\operatorname{d}\log|1-f_j|\wedge\operatorname{d}i\arg F_j(t_1)$$

$$-\log|g_j|\log\left|\frac{t_1-f_j}{t_1-1}\right|\operatorname{d}\log|1-f_j|\wedge\operatorname{d}i\arg F_j(t_2)$$

$$-\sum_{j,k,l}c_j s_{kj} t_{lj} P_{2,\text{Zag}}(f_j)\operatorname{d}i\arg\left(\frac{t_1-A_k}{t_1-1}\right)\wedge\operatorname{d}i\arg\left(\frac{t_2-A_l}{t_2-1}\right).$$

Observing that

$$\sum_j c_j\left(\log|g_j|\log|f_j|\operatorname{d}\log|1-f_j|\wedge\operatorname{d}i\arg\left(\frac{t_2-f_j}{t_2-1}\right)\right.$$

$$\left.-\log|g_j|\log|f_j|\operatorname{d}\log|1-f_j|\wedge\operatorname{d}i\arg F_j(t_2)\right)$$

$$=\sum_{j,k}c_j s_{kj}\log|g_j|\log|f_j|\operatorname{d}\log|1-f_j|\wedge\operatorname{d}i\arg\left(\frac{t_2-A_k}{t_2-1}\right)$$

we find for $\psi_1(0,t_2)$

$$\sum_j c_j\left(\frac{1}{6}\log|f_j|\sigma(g_j,1-f_j)\wedge\operatorname{d}\log\left|\frac{t_2-f_j}{t_2-1}\right|\right.$$

$$-\frac{1}{6}\log\left|\frac{t_2-f_j}{t_2-1}\right|\sigma(g_j,1-f_j)\wedge\operatorname{d}\log|f_j|$$

$$+\frac{1}{3}\log|g_j|\log|f_j|\operatorname{d}i\arg(1-f_j)\wedge\operatorname{d}\log\left|\frac{t_2-f_j}{t_2-1}\right|$$

$$-\frac{1}{3}\log|g_j|\log\left|\frac{t_2-f_j}{t_2-1}\right|\operatorname{d}i\arg(1-f_j)\wedge\operatorname{d}\log|f_j|$$

$$-\log|g_j|\log\left|\frac{t_2-f_j}{t_2-1}\right|\operatorname{d}\log|1-f_j|\wedge\operatorname{d}i\arg f_j$$

$$+\sum_k s_{kj}\log|g_j|\log|f_j|\operatorname{d}\log|1-f_j|\wedge\operatorname{d}i\arg\left(\frac{t_2-A_k}{t_2-1}\right)$$

$$\left.-\sum_k s_{kj} P_{2,\text{Zag}}(f_j)\operatorname{d}i\arg g_j\wedge\operatorname{d}i\arg\left(\frac{t_2-A_k}{t_2-1}\right)\right).$$

On the other hand, from writing down the explicit $\varepsilon_n$ for $n=1$, $2$ and $3$, we know that if we put $t_1 = 0$ in reg$^\sim$, we get dcor$_1$ with cor$_1(t_2)$ given by

$$\sum_j c_j \left( \varphi_1^{(3)} - \varphi_1^{(2)} \mathrm{d}i \arg F_j(t_2) \right) \wedge \mathrm{d}i \arg g_j$$

$$= \sum_j c_j \left( -P_{2,\mathrm{Zag}}(f_j) \mathrm{d}i \arg \left( \frac{t_2 - f_j}{t_2 - 1} \right) - \frac{2}{3} \log|1 - f_j| \log \left| \frac{t_2 - f_j}{t_2 - 1} \right| \mathrm{d} \log|f_j| \right.$$

$$\left. - \frac{1}{3} \log|1 - f_j| \log|f_j| \mathrm{d} \log \left| \frac{t_2 - f_j}{t_2 - 1} \right| + P_{2,\mathrm{Zag}}(f_j) \mathrm{d}i \arg F_j(t_2) \right) \wedge \mathrm{d}i \arg g_j$$

$$= \sum_j c_j \left( \sum_k s_{kj} - P_{2,\mathrm{Zag}}(f_j) \mathrm{d}i \arg \left( \frac{t_2 - A_k}{t_2 - 1} \right) - \frac{2}{3} \log|1 - f_j| \log \left| \frac{t_2 - f_j}{t_2 - 1} \right| \mathrm{d} \log|f_j| \right.$$

$$\left. - \frac{1}{3} \log|1 - f_j| \log|f_j| \mathrm{d} \log \left| \frac{t_2 - f_j}{t_2 - 1} \right| \right) \wedge \mathrm{d}i \arg g_j.$$

So $\psi_1(0, t_2) - \mathrm{cor}_1(t_2)$ equals

$$\sum_j c_j \left( \sum_k s_{kj} P_{2,\mathrm{Zag}}(f_j) \mathrm{d}i \arg \left( \frac{t_2 - A_k}{t_2 - 1} \right) + \frac{2}{3} \log|1 - f_j| \log \left| \frac{t_2 - f_j}{t_2 - 1} \right| \mathrm{d} \log|f_j| \right.$$

$$+ \frac{1}{3} \log|1 - f_j| \log|f_j| \mathrm{d} \log \left| \frac{t_2 - f_j}{t_2 - 1} \right| \right) \wedge \mathrm{d}i \arg g_j$$

$$\sum_j c_j \left( \frac{1}{6} \log|f_j| \sigma(g_j, 1 - f_j) \wedge \mathrm{d} \log \left| \frac{t_2 - f_j}{t_2 - 1} \right| \right.$$

$$- \frac{1}{6} \log \left| \frac{t_2 - f_j}{t_2 - 1} \right| \sigma(g_j, 1 - f_j) \wedge \mathrm{d} \log |f_j|$$

$$+ \frac{1}{3} \log|g_j| \log|f_j| \mathrm{d}i \arg(1 - f_j) \wedge \mathrm{d} \log \left| \frac{t_2 - f_j}{t_2 - 1} \right|$$

$$- \frac{1}{3} \log|g_j| \log \left| \frac{t_2 - f_j}{t_2 - 1} \right| \mathrm{d}i \arg(1 - f_j) \wedge \mathrm{d} \log|f_j|$$

$$- \log|g_j| \log \left| \frac{t_2 - f_j}{t_2 - 1} \right| \mathrm{d} \log|1 - f_j| \wedge \mathrm{d}i \arg f_j$$

$$+ \sum_k s_{kj} \log|g_j| \log|f_j| \mathrm{d} \log|1 - f_j| \wedge \mathrm{d}i \arg \left( \frac{t_2 - A_k}{t_2 - 1} \right)$$

$$\left. - \sum_k s_{kj} P_{2,\mathrm{Zag}}(f_j) \mathrm{d}i \arg g_j \wedge \mathrm{d}i \arg \left( \frac{t_2 - A_k}{t_2 - 1} \right) \right)$$

$$= \sum_{j,k} c_j t_{kj} \left( -2 P_{2,\mathrm{Zag}}(f_j) \mathrm{d}i \arg f_j + \log^2 |f_j| \mathrm{d} \log|1 - f_j| \right) \wedge \mathrm{d}i \arg \left( \frac{t_2 - A_k}{t_2 - 1} \right)$$

$$+ \sum_j c_j \left( \frac{1}{6} \log|f_j| \sigma(g_j, 1 - f_j) \wedge \mathrm{d} \log \left| \frac{t_2 - f_j}{t_2 - 1} \right| \right.$$

$$- \frac{1}{6} \log \left| \frac{t_2 - f_j}{t_2 - 1} \right| \sigma(g_j, 1 - f_j) \wedge \mathrm{d} \log|f_j|$$

$$+\frac{1}{3}\log|g_j|\log|f_j|\mathrm{d}i\arg(1-f_j)\wedge\mathrm{d}\log\left|\frac{t_2-f_j}{t_2-1}\right|$$

$$-\frac{1}{3}\log|g_j|\log\left|\frac{t_2-f_j}{t_2-1}\right|\mathrm{d}i\arg(1-f_j)\wedge\mathrm{d}\log|f_j|$$

$$-\log|g_j|\log\left|\frac{t_2-f_j}{t_2-1}\right|\mathrm{d}\log|1-f_j|\wedge\mathrm{d}i\arg f_j$$

$$+\frac{2}{3}\log|1-f_j|\log\left|\frac{t_2-f_j}{t_2-1}\right|\mathrm{d}\log|f_j|\wedge\mathrm{d}i\arg g_j$$

$$+\frac{1}{3}\log|1-f_j|\log|f_j|\mathrm{d}\log\left|\frac{t_2-f_j}{t_2-1}\right|\wedge\mathrm{d}i\arg g_j\Bigg)$$

where we used (4.10) and (4.11), as well as the fact that the map $[f]_2 \mapsto P_{2,\mathrm{Zag}}(f)$ factors through $M_{(2)}((F))$ (following (4.13)) and (4.13). Using integration by parts to get $\log\left|\frac{t_2-f_j}{t_2-1}\right|$ out, we find that this equals $\mathrm{d}\psi_2(t_2)$ with

$$\psi_2(t_2) = \sum_j c_j \left(\sum_k -2t_{kj}P_{3,\mathrm{Zag}}(f_j)\mathrm{d}i\arg\left(\frac{t_2-A_k}{t_2-1}\right)\right.$$

$$-\frac{1}{6}\log|f_j|\log\left|\frac{t_2-f_j}{t_2-1}\right|\sigma(g_j,1-f_j)$$

$$-\frac{1}{3}\log|g_j|\log|f_j|\log\left|\frac{t_2-f_j}{t_2-1}\right|\mathrm{d}i\arg(1-f_j)$$

$$+\frac{1}{3}\log|1-f_j|\log|f_j|\log\left|\frac{t_2-f_j}{t_2-1}\right|\mathrm{d}i\arg g_j\Bigg)$$

as all the other elements then cancel using only that $\pi_1$ of a holomorphic 2–form on $U$ vanishes identically.

We now have the (alternating) forms $\mathrm{reg}^\sim(t_1,t_2)$, $\psi_1(t_1,t_2)$, $\psi_2(t_2)$ and $\mathrm{cor}_1(t_2)$, all of which become identically zero if we put $t_2 = \infty$. If we let $\mathrm{cor}_2 = \sum_j c_j \varphi_2^{(3)}(f_j) \wedge \mathrm{d}i\arg g_j$ (which does not depend on $t_1$ or $t_2$), we have that they satisfy the relations

$$\mathrm{reg}^\sim(t_1,t_2) = \mathrm{d}\psi_1(t_1,t_2)$$
$$\mathrm{reg}^\sim(0,t_2) = \mathrm{d}\mathrm{cor}_1(t_2)$$
$$\mathrm{cor}_1(0) = \mathrm{d}\mathrm{cor}_2$$
$$\psi_1(0,t_2) - \mathrm{cor}_1(t_2) = \mathrm{d}\psi_2(t_2).$$

We now consider the following form, where $\rho_1(t)$ is a bump form around $t = 0$, $\rho_2(t_1,t_2)$ a bump form around the collection of $t_i = f_j$ or $t_i = 1$ (symmetric with respect to interchanging $t_1$ and $t_2$). It is easy to choose them in such a way that of $\rho_1(t_1)$, $\rho_1(t_2)$ and $\rho_2(t_1,t_2)$ at most two are non–zero at the same time. Let reg be the form (satifying (2.5)) given by

$$\begin{aligned}
&\mathrm{reg}^\sim - \mathrm{d}[\rho_1(t_1)\mathrm{cor}_1(t_2)] + \mathrm{d}[\rho_1(t_2)\mathrm{cor}_1(t_1)] \\
&+ \mathrm{d}\left[\rho_1(t_1)\mathrm{d}[\rho_1(t_2)\mathrm{cor}_2]\right] - \mathrm{d}\left[\rho_1(t_2)\mathrm{d}[\rho_1(t_1)\mathrm{cor}_2]\right] \\
&- \mathrm{d}\left[\rho_2(t_1,t_2)\psi_1(t_1,t_2)\right] \\
&+ \mathrm{d}\left[\rho_2(t_1,t_2)\left(\rho_1(t_1)\mathrm{cor}_1(t_2) - \rho_1(t_2)\mathrm{cor}_1(t_1)\right)\right] \\
&+ \mathrm{d}\left[\rho_2(t_1,t_2)\mathrm{d}[\rho_1(t_1)\psi_2(t_2)]\right] - \mathrm{d}\left[\rho_2(t_1,t_2)\mathrm{d}[\rho_1(t_2)\psi_2(t_1)]\right].
\end{aligned} \quad (4.15)$$

Here the first five terms are the original regulator as in (4.14), the remaining form the d of some 2–form vanishing at $t_i = 0, \infty$, lifting the regulator back to $X_V^2$, in fact to $(\mathbb{P}^1)_V^2$. (For checking that this is the case, note that the product $\rho_2(t_1, t_2)\rho_1(t_1)\rho_1(t_2)$ is identically zero by our choice of $\rho_2$, so that at most two of them are non–zero at any point of $X_U^2$.)

We now proceed to computing the integral in (4.8). Let $\Xi = \mathrm{d}i \arg t_1 \wedge \mathrm{d}i \arg t_2$. Some calculations using the formulae at the end of Section 2 give the following integrals.

$$\int_{X^2 \times S_x^1} \varepsilon_3(t_1, t_2, f) \wedge \mathrm{d}i \arg g \wedge \Xi = (2\pi i)^2 \int_{S_x^1} 2P_{3,\mathrm{Zag}}(f) \mathrm{d}i \arg g$$

$$\int_{X^2 \times S_x^1} -\log|g| \pi_3 \omega_3(t_1, t_2, f) = (2\pi i)^2 \int_{S_x^1} \log^2|f| \log|g| \mathrm{d}i \arg(1-f)$$

$$\int_{X^2 \times S_x^1} \mathrm{d}i \arg F(t_1) \wedge \varepsilon_2(t_2, f) \wedge \mathrm{d}i \arg g \wedge \Xi = 0$$

$$\int_{X^2 \times S_x^1} \mathrm{d}i \arg F(t_1) \wedge \log|g| \pi_2 \omega_2(t_2, f) \wedge \Xi = 0$$

$$\int_{X^2 \times S_x^1} \mathrm{d}i \arg F(t_1) \wedge \mathrm{d}i \arg F(t_2) \wedge \sigma(1-f, g) \wedge \Xi = 0$$

$$\int_{X^2 \times S_x^1} -\mathrm{d}[\rho_2(t_1, t_2)\psi_1(t_1, t_2)] \wedge \Xi = -4\pi i \int_{X^1 \times S_x^1} \rho_2(0, t)\psi_1(0, t) \wedge \mathrm{d}i \arg t$$

$$\int_{X^2 \times S_x^1} \mathrm{d}[\rho_2(t_1, t_2)(\rho_1(t_1)\mathrm{cor}_1(t_2) - \rho_1(t_2)\mathrm{cor}_1(t_1))] \wedge \Xi =$$

$$4\pi i \int_{X^1 \times S_x^1} \rho_2(0, t) \mathrm{cor}_1(t) \wedge \mathrm{d}i \arg t$$

$$\int_{X^2 \times S_x^1} \mathrm{d}[\rho_2(t_1, t_2)\mathrm{d}[\rho_1(t_1)\psi_2(t_2)]] - \mathrm{d}[\rho_2(t_1, t_2)\mathrm{d}[\rho_1(t_2)\psi_2(t_1)]] \wedge \Xi =$$

$$4\pi i \int_{X^1 \times S_x^1} \rho_2(0, t) (\psi_1(0, t) - \mathrm{cor}_1(t)) \wedge \mathrm{d}i \arg t + (2\pi i)^2 \int_{S_x^1} 4P_{2,\mathrm{Zag}}(f) \mathrm{d}i \arg g$$

(For those computations, note that $S_x^1$ has dimension one, so that a lot of the contributions actually vanish identically on $X_{S_x^1}^2$.) The first five lines here suffice to compute the contribution of the first two lines in (4.15), as those equal (4.14).

Putting everything together we find that (4.8) equals

$$\frac{-1}{2\pi i} \int_{S_x^1} \sum_j c_j \left( 6P_{3,\mathrm{Zag}}(f_j) \mathrm{d}i \arg g_j + \log^2|f_j| \log|g_j| \mathrm{d}i \arg(1-f_j) \right). \qquad (4.16)$$

We now rewrite the form in (4.16) in order to compute its residues. Using (2.4) and (4.12), we find

$$\sum_j c_j \left( 6P_3^{\mathrm{mod}}(f_j) \mathrm{d}i \arg g_j + \log^2|f_j| \sigma(g_j, 1-f_j) \right) =$$

$$\sum_j c_j \left( 6P_3^{\mathrm{mod}}(f_j) \mathrm{d}i \arg g_j + \log|f_j| \log|g_j| \sigma(f_j, 1-f_j) \right).$$

Subtracting

$$\sum_j \mathrm{d}c_j \left(\log|f_j|\log|g_j|P_{2,\mathrm{Zag}}(f_j)\right) = \sum_j c_j \left(2\log|g_j|P_{2,\mathrm{Zag}}(f_j)\mathrm{d}\log|f_j|\right.$$
$$\left. + \log|f_j|\log|g_j|\sigma(f_j, 1-f_j)\right)$$

this transforms into

$$\sum_j c_j \left(6 P_3^{\mathrm{mod}}(f_j)\mathrm{d}i\arg g_j - 2\log|g_j|P_2^{\mathrm{mod}}(f_j)\mathrm{d}\log|f_j|\right). \tag{4.17}$$

The integral in (4.16) therefore yields

$$-6 \sum_j c_j \mathrm{ord}_x(g_j) P_3^{\mathrm{mod}}(f_j(x)).$$

This is the regulator of the element

$$\pm 3 \sum_j c_j \delta_x([f_j]_3 \otimes g_j),$$

see Theorem 2.3. Because the regulator is injective on $K_5^{(3)}(k(x))$, this proves that the following diagram commutes (up to sign and up to $\partial\left(K_5^{(3)}(k) \cup F_{\mathbb{Q}}^*\right)$ in the lower right hand corner):

$$\begin{array}{ccc}
H^2(\mathcal{M}_{(4)}^\bullet(F)) & \longrightarrow & K_6^{(4)}(F)/K_4^{(2)}(F) \cup K_2^{(2)}(F) \\
\downarrow{\scriptstyle 3\delta} & & \downarrow{\scriptstyle \partial} \\
\coprod_{x \in C^{(1)}} H^1(\widetilde{\mathcal{M}}_{(3)}^\bullet(k(x))) & \xrightarrow{\sim} & \coprod_{x \in C^{(1)}} K_5^{(3)}(k(x))
\end{array}$$

**Remark 4.20** Because the form $\psi$ appearing in (4.17) (or (4.16)) has the same residue (modulo the residue of the regulator of $K_5^{(3)}(k) \cup F_{\mathbb{Q}}^*$) as the regulator of $\alpha$ in the localization sequence

$$0 \to H^1_{\mathrm{dR}}(C; \mathbb{R}(3))^+ \to H^1_{\mathrm{dR}}(F; \mathbb{R}(3))^+ \to \coprod H^0_{\mathrm{dR}}(k(x); \mathbb{R}(2))^+$$

they differ by an element in $H^1_{\mathrm{dR}}(C; \mathbb{R}(3))^+ + \mathrm{reg}(K_5^{(3)}(k) \cup F_{\mathbb{Q}}^*)$. Using integration by parts and Stokes' theorem it is not hard to check that

$$\int_C \psi \wedge \overline{\omega} = \int_C \mathrm{reg}(\alpha) \wedge \overline{\omega}$$

with the last given by Theorem 3.4. By Remark 3.1 $\psi$ is an explicit representative of the regulator of $\sum_j c_j[f_j]_3 \otimes g_j$ in $H^1_{\mathrm{dR}}(C; \mathbb{R}(3))$, modulo the regulator of $K_5^{(3)}(k) \cup F_{\mathbb{Q}}^*$.

# 5 Connections with Goncharov's work

In this section, we start with showing how the work in the previous two sections, together with the work of Goncharov ([G1] and [G3], see [G2, Section 8] for an overview of the results without proofs) leads to a complete description of the image of the regulator map on $K_4^{(3)}(C)$ and $K_6^{(4)}(C)$. In particular, this proves a conjecture of Goncharov for those cases ([G3, Conjecture 1.5] or [G1, Conjecture 1.6]). We also sketch how, assuming some conjectures, the relation with results as in Goncharov's work for higher $K$–groups would work out.

In [G2, § 6], Goncharov defined the following complexes $\Gamma(F, n)$ (in degree $1, \ldots, n$), given by

$$B_n(F) \to B_{n-1}(F) \otimes F_{\mathbb{Q}}^* \to \ldots \to B_2(F) \otimes \bigwedge^{n-2} F_{\mathbb{Q}}^* \to \bigwedge^n F_{\mathbb{Q}}^*$$

and, for each $x \in C^{(1)}$, $\Gamma(k(x), n-1))$ (in degrees $1, \ldots, n-1$), given by

$$B_{n-1}(k(x)) \to \ldots \to B_2(k(x)) \otimes \bigwedge^{n-3} k(x)_{\mathbb{Q}}^* \to \bigwedge^{n-1} k(x)_{\mathbb{Q}}^*$$

Here for any infinite field $F$, $B_k(F)$ is a $\mathbb{Q}$–vector space generated by elements $\{f\}_k$ with $f \in F \cup \{\infty\}$, modulo certain (inductively defined) relations. All maps are given by

$$\{f\}_k \otimes g_1 \wedge \ldots \wedge g_{n-k} \mapsto \{f\}_{k-1} \otimes f \wedge g_1 \wedge \ldots \wedge g_{n-k}.$$

There is a map

$$\Gamma(F, n) \to \coprod_{x \in C^{(1)}} \Gamma(k(x), n-1)[-1] \tag{5.1}$$

given by

$$\{f\}_k \otimes g_1 \wedge \ldots \wedge g_{n-k} \mapsto \{f(x)\}_k \otimes \partial_{n-k,x}(g_1 \wedge \ldots \wedge g_{n-k})$$

with $\partial_{m,x}$ the unique map $\bigwedge^m F_{\mathbb{Q}}^* \to \bigwedge^{m-1} k(x)_{\mathbb{Q}}^*$ determined as in Proposition 4.1

$$\pi_x \wedge u_1 \wedge \ldots \wedge u_{k-1} \mapsto u_1(x) \wedge \cdots \wedge u_{k-1}(x)$$
$$u_1 \wedge \ldots \wedge u_k \mapsto 0$$

if all $u_i$ are units at $x$ and $\pi_x$ is a uniformizer at $x$. $\Gamma(C, n)$ is defined as the mapping cone of (5.1). Goncharov also defines complexes $\Gamma'(F, n)$, $\Gamma'(k(x), n-1)$ for $n = 3$ and 4, and constructs maps as in (5.1). The complexes $\Gamma'$ have the same shape as the complexes $\Gamma$ with the same maps between them, but the $B_k(F)$ get replaced with $B'_k(F)$, generated by $F \cup \{\infty\}$, but with *explicit* relations between the generators. $\Gamma'(C, n)$ is defined as the mapping cone, defined by the corresponding $\Gamma'$ complexes in (5.1). Goncharov also constructs a map

$$K_{2n}(C) \to H^2(\Gamma'(C, n+1)) \tag{5.2}$$

for $n = 2$ or 3, and shows that the Beilinson regulator factors through this map. We summarize part of his results in a form suitable for our needs.

**Theorem 5.1 (Goncharov)** Let $\omega$ be a global holomorphic 1–form on $C$. Then for $n = 2$ or 3, the regulator map

$$K_{2n}^{(n+1)}(C) \to H_{\mathrm{dR}}^1(C; \mathbb{R}(n))^+$$

can be extended over $K_{2n}^{(n+1)}(C) \to H^2(\Gamma'(C, n+1))$ to $H^2(\Gamma'(C, n+1)) \to H^1_{\mathrm{dR}}(C; \mathbb{R}(n))^+$.
For $\omega$ a holomorphic 1–form, the composition

$$H^2(\Gamma'(C, n+1)) \longrightarrow H^1_{\mathrm{dR}}(C; \mathbb{R}(n))^+ \xrightarrow{\int_C \cdots \wedge \overline{\omega}} \mathbb{R}(1)$$

is given by mapping $\{f\}_n \otimes g$ to

$$c_n \int_C \log|g| \log^{n-2}|f| \left(\log|1-f| \mathrm{d}\log|f| - \log|f| \mathrm{d}\log|1-f|\right) \wedge \overline{\omega}$$

for some non–zero rational constant $c_n$.

For $n = 2$, this is proved in [G3]. There Theorem 2.2 constructs the map (5.2), the extension of the regulator is given just before Theorem 3.1, which states that the extension of the regulator coincides with Beilinson's regulator on $K_4^{(3)}(C)$. Finally, Theorem 3.3 gives the formula for the regulator integral. For $n = 3$, the corresponding results can be found in [G1], namely Theorems 4.2, 5.3 and 5.5.

**Lemma 5.2** There is a map
$$B'_2(F) \to \widetilde{M}_{(2)}(F)$$
given by sending $\{x\}_2$ to $[x]_2$.

**Proof** $B'_2(F)$ is a free $\mathbb{Q}$–vector space on elements $\{x\}_2$ with $x \in F^* \setminus \{1\}$, modulo the relations
$$\{x\}_2 + \{y\}_2 + \{\frac{1-x}{1-xy}\}_2 + \{1-xy\}_2 + \{\frac{1-y}{1-xy}\}_2 = 0$$

It is known (see [Su, Lemmas 1.2 and 1.4]) that one then also has the following relations
$$\{x\}_2 + \{1-x\}_2 = 0 \qquad \text{and} \qquad \{x\}_2 + \{1/x\}_2 = 0$$

We have to show that the corresponding relations hold in $\widetilde{M}_{(2)}(F)$. We start with the last two. The relation $[x]_2 + [1/x]_2 = 0$ holds in $\widetilde{M}_{(2)}(F)$ by definition. The element $[x]_2 + [1-x]_2$ lies in $H^1(\widetilde{\mathcal{M}}^\bullet_{(2)}(F))$ and is a pull back from an element in $H^1(\widetilde{\mathcal{M}}^\bullet_{(2)}(\mathbb{Q}(t)))$, which injects into $K_3^{(2)}(\mathbb{Q}(t))$ by Theorem 2.3. But $K_3^{(2)}(\mathbb{Q}(t)) = K_3^{(2)}(\mathbb{Q}) = 0$ so this element is zero. For the actual five–term relation, observe that modulo the last two relations (for $[\cdots]_2$ instead of $\{\cdots\}_2$), the first corresponds to the relation in $M_{(2)}(F)$ given by

$$-[x^{-1}]_2 - [1-y]_2 + [\frac{1-xy}{x}]_2 + [\frac{x}{1-xy}]_2 - [1 - \frac{1-x}{1-xy}]_2 - [\frac{1}{1-xy}]_2 + [\frac{1-y}{1-xy}]_2.$$

The construction of the complexes as sketched in Section 2 gives that the lift of those elements are given by
$$\sum_j [f_j]_2 + \sum_j (1-f_j) \cup F_j(t)$$

where $F_j$ is the function expressing $f_j$ in elements of a chosen basis in $\{f_j\}$. In order to show that those equal zero, we work universally, i.e., we work over the base scheme

$$Z = \mathrm{Spec}\,(\mathbb{Q}[X, Y, (1-X)^{-1}, (1-Y)^{-1}, (1-XY)^{-1}])$$

and we want to show that we are pulling back a universal element, obtained by replacing $x$ with $X$ and $y$ with $Y$, from $K_3^{(2)}(Z) \cong K_3^{(2)}(\mathbb{Q}) = 0$ via the map $x \mapsto X$, $y \mapsto Y$. We can pull back directly where all $f_j \neq 1$, i.e., pull back from the open part $Z'$ of $Z$ where $1-X-XY \neq 0$. But $K_3^{(2)}(Z') = 0$ as well. If $1-x-xy = 0$, using the relations $\{x\}_2 + \{1-x\}_2$ and $\{x\}_2 + \{1/x\}_2$ which we know already, the relation reduces to to $\{x^2\}_2 = 2\{x\}_2 + 2\{-x\}_2$. One proves this one in a similar way over Spec $(\mathbb{Q}[X, (1-X^2)^{-1}])$.

**Remark 5.3** In fact $B_2(F) = B_2'(F)$. Namely, let $F$ be any infinite field, and suppose $\alpha \in \operatorname{Ker}\left(\mathrm{d}: \mathbb{Q}[F(T) \cup \{\infty\}] \to \bigwedge^2 F(t)_{\mathbb{Q}}^*\right)$. By Suslin's work, this yields an element in $K_3^{(2)}(F(t)) \cong K_3^{(2)}(F)$. But modulo the five term relations, we can rewrite this to $\alpha \equiv \beta$ with $\beta \in K_3^{(2)}(F) \cong \operatorname{Ker}\left(\mathrm{d}: B_2(F) \to \bigwedge^2 F_{\mathbb{Q}}^*\right)$. Then $\alpha(0) \equiv \beta$ modulo the relations, as one checks by a case by case check depending on the zeroes and poles of the functions involved. Of course this works for $\alpha(1)$ as well, so in total $\alpha(0) - \alpha(1)$ is in the (degenerate) relations in $B_2$, hence is zero in $B_2(F)$.

We use Lemma 5.2 to link our results with Theorem 5.1, beginning with the case $n = 2$.

**Theorem 5.4** The maps in (5.2), Lemma 5.2 and $\varphi_{(3)}^2$ give maps

$$K_4^{(3)}(C) \to H^2(\Gamma'(C, 3)) \to H^2(\widetilde{\mathcal{M}}_{(3)}^\bullet(C)) \to K_4^{(3)}(C) + K_3^{(2)}(k) \cup F_{\mathbb{Q}}^*.$$

Viewing this last group as inside $K_4^{(3)}(F)$, the composition of those maps with the regulator integral associated to $\omega$ is given by Theorem 3.4 on $H^2(\widetilde{\mathcal{M}}_{(3)}^\bullet(C))$ and by Theorem 5.1 on $H^2(\Gamma'(C, 3))$. In particular, all those groups, as well as the group $H^2(\mathcal{M}_{(3)}^\bullet(C)) \cong H^2(\widetilde{\mathcal{M}}_{(3)}^\bullet(C))$, have the same image in $\mathbb{R}(1)$ under the regulator integral.

**Proof** That we get the map from $H^2(\Gamma'(C, 3))$ to $H^2(\widetilde{\mathcal{M}}_{(3)}^\bullet(C))$ is clear from Lemma 5.2. The form of the regulator integral on $H^2(\widetilde{\mathcal{M}}_{(3)}^\bullet(C))$ was stated in Theorem 3.4, and that the composition of this with the map $H^2(\Gamma'(C, 3)) \to H^2(\widetilde{\mathcal{M}}_{(3)}^\bullet(C))$ coincides with the formulas in Theorem 5.1 up to a non–zero rational number is clear. Then, fixing $\omega$, we get that the images in $\mathbb{R}(1)$ have the relations:

$$\operatorname{Image}(K_4^{(3)}(C)) \subset \operatorname{Image}(H^2(\Gamma'(C, 3))) \subset \operatorname{Image}(H^2(\widetilde{\mathcal{M}}_{(3)}^\bullet(C)))$$

$$\subset \operatorname{Image}(K_4^{(3)}(C) + K_3^{(2)}(k) \cup F_{\mathbb{Q}}^*) = \operatorname{Image}(K_4^{(3)}(C))$$

because the regulator integral vanishes on $K_3^{(2)}(k) \cup F_{\mathbb{Q}}^*$ by Proposition 3.3.

**Corollary 5.5** The groups $K_4^{(3)}(C)$, $H^2(\Gamma'(C, 3))$ and, in case $K_3^{(2)}(k) = 0$, the group $H^2(\mathcal{M}_{(3)}^\bullet(C)) \cong H^2(\widetilde{\mathcal{M}}_{(3)}^\bullet(C))$ have the same image in $H_{\mathrm{dR}}^1(C; \mathbb{R}(2))^+$ under the regulator map. The same holds true without assuming that $K_3^{(2)}(k) = 0$ if we use the modified version of $\varphi_{(2)}^3$ as described in Corollary 4.5.

**Proof** This is clear from Theorem 5.4, as the regulator integrals form the dual space of $H_{\mathrm{dR}}^1(C; \mathbb{R}(2))^+$, see Remark 3.1, and the difference between using $\varphi_{(3)}^2$ and its modification lies in $K_3^{(2)}(k) \cup F_{\mathbb{Q}}^*$, on which the regulator integral vanishes by Proposition 3.3.

We now turn towards $n = 3$. Because the natural map $H^2(\mathcal{M}^\bullet_{(4)}(F)) \to H^2(\widetilde{\mathcal{M}}^\bullet_{(4)}(F))$ is a surjection as described in Section 2, we get a surjection $H^2(\mathcal{M}^\bullet_{(4)}(C)) \to H^2(\widetilde{\mathcal{M}}^\bullet_{(4)}(C))$. In particular, those two groups have the same image under the regulator integral as in Theorem 3.4.

We recall the definition of the group $B'_3(k)$: it is the free $\mathbb{Q}$–vector space with generators $\{x\}_3$ for $x \in F^* \setminus \{1\}$, and relations

$$\sum_{i=1}^{3} \left( \{\alpha_i\}_3 + \{\beta_i\}_3 - \{\frac{\beta_i}{\alpha_{i-1}}\}_3 + \{\frac{\beta_i}{\alpha_{i-1}\alpha_i}\}_3 + \{\frac{\alpha_i\beta_{i-1}}{\beta_{i+1}}\}_3 \right. \tag{5.3}$$
$$\left. + \{-\frac{\beta_i}{\alpha_i\beta_{i-1}}\}_3 - \{\frac{\alpha_i\alpha_{i-1}\beta_{i+1}}{\beta_i}\}_3 \right) - 3\{1\}_3 + \{-\alpha_1\alpha_2\alpha_3\}_3.$$

Here $\beta_i = 1 - \alpha_i(1 - \alpha_{i-1})$ with indices taken modulo 3.

**Theorem 5.6** Let $C$ be a smooth, proper, geometrically irreducible curve over the number field $k$, with function field $F$. Then the groups

$$K_6^{(4)}(C), \ H^2(\Gamma'(C,4)), \ H^2(\mathcal{M}^\bullet_{(3)}(C)) \text{ and } H^2(\widetilde{\mathcal{M}}^\bullet_{(3)}(C))$$

all have the same image under the regulator integral, given by Theorem 3.4.

**Proof** By Corollary 4.8, $H^2(\mathcal{M}^\bullet_{(4)}(C))$ maps to $K_6^{(4)}(C) + K_5^{(3)}(k) \cup F_\mathbb{Q}^* \subset K_6^{(4)}(F)$. The regulator integral is zero on $K_5^{(3)}(k) \cup F_\mathbb{Q}^*$ by Proposition 3.3, whence

$$\text{Image}(H^2(\widetilde{\mathcal{M}}^\bullet_{(4)}(C))) = \text{Image}(H^2(\mathcal{M}^\bullet_{(4)}(C))) \subset \text{Image}(K_6^{(4)}(C))$$

in $\mathbb{R}(1)$.

From Goncharov's work as quoted in Theorem 5.1 we get a map

$$K_6^{(4)}(C) \to H^2(\Gamma'(C,4)).$$

By Lemma 5.2 we have a map
$$B'_2(F) \to \widetilde{M}_{(2)}(F).$$

In [dJ1, p. 241] a map
$$B'_3(k(x)) \to \widetilde{M}_{(3)}(k(x))$$

was created by mapping $\{y\}_3$ to $[y]_3$. This is well defined because d is zero on the relations in (5.3) so they give rise to an element in $H^1(\widetilde{\mathcal{M}}^\bullet_{(3)}(k(x))) \subset K_5^{(3)}(k(x))$. Because the regulator for the embedding $\sigma$ of $k(x)$ into $\mathbb{C}$ is given by mapping $[y]_3$ to a non–zero multiple of $P_3^{\text{mod}}(\sigma(y))$ the function $P_3^{\text{mod}}$ vanishes on the elements in (5.3) and the regulator is injective, the elements in (5.3) go to zero in $\widetilde{M}_{(3)}(k(x))$ and our map is well–defined. Using those two maps we see that if $\sum_j c_j\{f_j\}_3 \otimes g_j$ is an element of $H^2(\Gamma'(C,4))$, then $\sum_j c_j[f_j]_3 \otimes g_j$ is an element of $H^2(\widetilde{\mathcal{M}}^\bullet_{(4)}(C))$. So we get the inclusion of images under the regulator integrals

$$\text{Image}(K_6^{(4)}(C)) \subset \text{Image}(H^2(\Gamma'(C,4))) \subset \text{Image}(H^2(\widetilde{\mathcal{M}}^\bullet_{(4)}(C))) \subset \text{Image}(K_6^{(4)}(C))$$

**Remark 5.7** If we use the lifted version $\varphi_{(4)}^2 : H^2(\mathcal{M}_{(4)}^\bullet(C)) \to K_6^{(4)}(C)/K_4^{(2)}(F) \cup K_2^{(2)}(F)$ as in Corollary 4.9, and combine this with the regulator map to $H^1_{\mathrm{dR}}(C; R(3))^+$, the resulting map factors through the projection $H^2(\mathcal{M}_{(4)}^\bullet(C)) \to H^2(\widetilde{\mathcal{M}}_{(4)}^\bullet(C))$. To see this, note that the kernel of this projection consists of elements $\alpha$ of the form $\sum_j n_j([f_j]_3 - [1/f_j]_3) \otimes g_j$, with $\sum_j n_j([f_j]_2 + [1/f_j]_2) \otimes f_j \wedge g_j = 0$ in $M_{(2)} \otimes \bigwedge^2 F_\mathbb{Q}^*$. By Remark 3.1, it is enough to check that the regulator integrals all vanish on the regulator of such $\alpha$, but this is part of Theorem 3.4.

Using the lifted versions of $\varphi_{(4)}^2$ and Remark 5.7, we get the following Corollary.

**Corollary 5.8** The groups $K_6^{(4)}(C)$, $H^2(\Gamma'(C, 4))$, $H^2(\mathcal{M}_{(4)}^\bullet(C))$ and $H^2(\widetilde{\mathcal{M}}_{(4)}^\bullet(C))$ have the same image in $H^1_{\mathrm{dR}}(C; \mathbb{R}(3))^+$ under the regulator map.

**Proof** This is immediate from Theorem 5.6 because the regulator integrals are dual to $H^1_{\mathrm{dR}}(C; \mathbb{R}(3))^+$, see Remark 3.1.

**Remark 5.9** According to the Beilinson conjectures, for $n \geq 2$, the regulator map $K_{2n}^{(n+1)}(C) \to H^1_{\mathrm{dR}}(C; \mathbb{R}(n))^+$ should be an injection, and in fact an isomorphism after tensoring the left hand side with $\mathbb{R}$. If that is the case, then Theorem 5.4 and Theorem 5.6 give, in principle, a complete combinatorial description in terms of generators of $K_4^{(3)}(C)$ resp. $K_6^{(4)}(C)$, and some (but not necessarily all) relations (i.e., $\{f\}_2 \otimes f$ and $\{f\}_3 \otimes f$ are zero).

**Remark 5.10** One could try to check the explicit relations of $B_3'(F)$ in $\widetilde{M}_{(3)}(F)$ along the lines of Lemma 5.2, in order to get a map from $H^2(\Gamma'(C, 4)) \to H^2(\widetilde{\mathcal{M}}_{(4)}^\bullet(C)))$. Due to the size of the relations involved, the author has not tried to do this. Note also that that would still not give as a map from $K_6^{(4)}(C) \to K_6^{(4)}(C) + K_5^{(3)}(k) \cup F_\mathbb{Q}^*$ similar to the maps in Theorem 5.4, as the map $\varphi_{(4)}^2$ from $H^2(\widetilde{\mathcal{M}}_{(4)}^\bullet(C))$ to $K_6^{(4)}(C) + K_5^{(3)}(k) \cup F_\mathbb{Q}^*$ depends on the Beilinson–Soulé conjecture as explained in Section 2. Thus the results for $n = 3$ are necessarily weaker than those for $n = 2$.

**Remark 5.11** One can give a more general proof of the existence of a map $B_n(F) \to \widetilde{M}_{(n)}(F)$ for all $n \geq 2$, but it becomes dependent on conjectures. Namely, assume that

1) $F$ is the function field of a smooth, projective, geometrically irreducible variety $Z$ over the number field $k$;
2) $B_n(F)$ is a quotient of the free $\mathbb{Q}$–vector space on elements $\{x\}_n$ with $x \in F^* \setminus \{1\}$, with elements of the form $\sum_j c_j\{f_j(x_1, \ldots, x_m)\}_n = 0$ for rational numbers $c_j$, and rational functions $f_j$ on $Z$ with coefficients in a number field $k$. Assume moreover that there exists a Zariski open part $U$ of $Z$ such that for all $y$ closed in $U$, the function $\sum_j c_j P_n^{\mathrm{mod}}(\sigma(f_j(y)))$ vanishes identically for all embeddings of $k(y)$ into $\mathbb{C}$;
3) the Beilinson–Soulé conjecture is true for general fields of characteristic zero: $K_n^{(p)}(F) = 0$ if $2p - n \leq 0$ and $n > 0$;
4) For a smooth, geometrically irreducible variety $Z$ over a number field $k$, $K_{2n-1}^{(n)}(Z) \cong K_{2n-1}^{(n)}(k)$ by pullback from the base (which is part of the Beilinson conjectures).

Then proceeding by induction, assume that we have defined a map $B_{n-1}(F) \to \widetilde{M}_{(n-1)}(F)$ by $\{x\}_{n-1} \mapsto [x]_{n-1}$, so that the diagram

$$\begin{array}{ccc} <\{x\}_n, x \in F^* \setminus \{1\}> & \longrightarrow & B_{n-1}(F) \otimes F_{\mathbb{Q}}^* \\ \downarrow & & \downarrow \\ \widetilde{M}_{(n)}(F) & \longrightarrow & \widetilde{M}_{(n-1)}(F) \otimes F_{\mathbb{Q}}^* \end{array}$$

(resp.

$$\begin{array}{ccc} <\{x\}_n, x \in F^* \setminus \{1\}> & \longrightarrow & \bigwedge^2 F_{\mathbb{Q}}^* \\ \downarrow & & \downarrow \\ \widetilde{M}_{(2)}(F) & \longrightarrow & \bigwedge^2 F_{\mathbb{Q}}^* \end{array}$$

for $n = 2$) commutes. Then we have to check that for any relation $\sum c_j \{f_j(x_1, \ldots, x_m)\}_n = 0$, the corresponding relation $\sum c_j [f_j(x_1, \ldots, x_m)]_n = 0$ holds in $\widetilde{M}_{(n)}(F)$. The element $\sum c_j [f_j(x_1, \ldots, x_m)]_n$ defines an element $\alpha$ in $H^1(\widetilde{\mathcal{M}}_{(n)}^\bullet(F))$, injecting into $K_{2n-1}^{(n)}(F)$. Using the spectral sequence

$$E_1^{p,q} = \coprod_{x \in X^{(p)}} K_{-p-q}^{(n-p)}(k(x)) \Rightarrow K_{-p-q}^{(n)}(Z)$$

(see [So, Théorème 4 (iii)]) the Beilinson-Soulé conjecture implies that then $K_{2n-1}^{(n)}(F) \cong K_{2n-1}^{(n)}(Z)$, and the Beilinson conjectures imply that $K_{2n-1}^{(n)}(k)$ by pullback from the base. The remarks in Section 2 show that $\varphi_{(n)}^1$ is in fact defined over some Zariski open subset $U$ of $Z$, and we have $K_{2n-1}^{(n)}(F) \cong K_{2n-1}^{(n)}(U) \cong K_{2n-1}^{(n)}(k)$ as well. We can select a point $y$ in $U$ such that $\varphi_{(n)}^1(\alpha)$ can be pulled back to $y$, mapping $\alpha$ to an element in $K_{2n-1}^{(n)}(k(y))$, namely the image from the corresponding element in $K_{2n-1}^{(n)}(k) \cong K_{2n-1}^{(n)}(U)$. Because the map $K_{2n-1}^{(n)}(k) \to K_{2n-1}^{(n)}(k(y))$ is injective, we can check that the image of $\alpha$ is zero by computing the regulator map, which according to Theorem 2.3 is given by computing $\sum_j c_j P_n^{\mathrm{mod}}(\sigma(f_j(y)))$ for all embeddings $\sigma$ of $k(y)$ into $\mathbb{C}$. This vanishes by our assumptions.

### References

[De]   C. Deninger. Higher order operations in Deligne cohomology. *Inventiones Mathematicae*, 120: 289—315, 1995.
[G1]   A.B. Goncharov. Special values of Hasse-Weil $L$–functions and generalized Eisenstein-Kronecker series. Preprint, 1994.
[G2]   A.B. Goncharov. Polylogarithms in Arithmetic and Geometry. In: Proc. International congress of Mathematicians (Zürich), 1994, pp. 374—387.
[G3]   A.B. Goncharov. Deninger's conjecture on $L$-functions of elliptic curves at $s = 3$. Preprint, 1994. Preprint MPI/96-6 (January 96), 38 pages. To appear in the special volume dedicated to Manin's 60-th birthday (1997).
[G–H] P. Griffiths and J. Harris.Principles of Algebraic Geometry. New York: John Wiley and sons, 1978.

[dJ1]   R. de Jeu. Zagier's Conjecture and Wedge Complexes in Algebraic $K$–theory. *Compositio Mathematica*, 96:197—247, 1995.
[dJ2]   R. de Jeu. On $K_4^{(3)}$ of curves over number fields. *Inventiones Mathematicae*, 125: 523—556, 1996.
[So]    C. Soulé. Opérations en $K$–théorie algébrique. *Canadian Journal of Mathematics*, XXXVII No.3: 488—550, 1985.
[Su]    A. A. Suslin. $K_3$ of a field and the Bloch group. *Trudy Math. Inst. Steklov*, 189:180—199, 1990. Russian, English translation *Proc. Steklov Inst. Math.*, vol. 183/4, pages 217—239, 1991.
[Za]    D. Zagier. Polylogarithms, Dedekind Zeta Functions, and the Algebraic $K$–Theory of Fields. In: G. van der Geer, F. Oort, J. Steenbrink, (eds.) Arithmetic Algebraic Geometry. Basel: Birkhäuser, 1991, pp. 391—430